\documentstyle[12pt,amssymb,amscd]{amsart}  

\newtheorem{thm}{Theorem}[section]
\newtheorem{cor}[thm]{Corollary}
\newtheorem{prop}[thm]{Proposition}
\newtheorem{lemma}[thm]{Lemma}
\newtheorem{claim}[thm]{Claim}
\newtheorem{conj}[thm]{Conjecture}

\theoremstyle{remark}
\newtheorem{remark}[thm]{Remark}
\newtheorem{example}[thm]{Example}

\theoremstyle{definition}
\newtheorem{defn}[thm]{Definition}

\numberwithin{equation}{section}

\newcommand{\bbL}{{\Bbb L}}
\newcommand{\bbR}{{\Bbb R}}

\newcommand{\bbZ}{{\Bbb Z}}
\newcommand{\bbD}{{\Bbb D}}

\newcommand{\cF}{{\cal F}}
\newcommand{\cG}{{\cal G}}

\newcommand{\cL}{{\cal L}}
\newcommand{\cM}{{\cal M}}
\newcommand{\cN}{{\cal N}}

\newcommand{\cV}{{\cal V}}
\newcommand{\cA}{{\cal A}}

\newcommand{\cI}{{\cal I}}
\newcommand{\cC}{{\cal C}}

\newcommand{\cK}{{\cal K}}
\newcommand{\cH}{{\cal H}}

\newcommand{\Sh}{\operatorname{Sh}}
\newcommand{\supp}{\operatorname{Supp}}

\newcommand{\Ker}{\operatorname{Ker}}

\newcommand{\R}{\operatorname{\bold R}}
\newcommand{\LL}{\operatorname{\bold L}}
\newcommand{\Hom}{\operatorname{Hom}}
\newcommand{\Ext}{\operatorname{Ext}}
\newcommand{\Tor}{\operatorname{Tor}}
\newcommand{\isomo}{\overset{\sim}{=}}

\newcommand{\Id}{\operatorname{Id}}

\newcommand{\Star}{\operatorname{Star}}
\newcommand{\Span}{\operatorname{Span}}

\newcommand{\Or}{{\underline{o}}}

\thanks
{This research was supported in part by the NSF}

\begin{document}

\title{Intersection cohomology on nonrational polytopes}

\author{Paul Bressler}
\address{I.H.E.S., 35, route de Chartres, 91440 Bures-sur-Yvette, France}
\email{bressler@@ihes.fr}
\author{Valery A.~Lunts}
\address{Department of Mathematics, Indiana University,
Bloomington, IN 47405, USA}
\email{vlunts@@indiana.edu}

\maketitle

\section{Introduction}
For an $n$-dimensional convex polytope $Q$ R.~Stanley ([S])
 defined a set of integers
$h(Q)= (h_0(Q), h_1(Q),\dots, h_n(Q))$ - the ``generalized $h$-vector'' - 
which are supposed to be the intersection cohomology Betti numbers of the
toric variety $X_Q$ corresponding to $Q$. In case $Q\subset\bbR^n$
is a rational polytope the variety $X_Q$ indeed exists, and it is known ([S])
that  $h_i(Q)=\dim IH^{2i}(X_Q)$. Thus, for a {\em rational polytope $Q$},
the integers $h_i(Q)$ satisfy 
\begin{enumerate}
\item $h_i(Q)\geq 0$,
\item $h_i(Q)=h_{n-i}(Q)$ (Poincar\'e duality),
\item $h_0(Q)\leq h_1(Q)\leq ... \leq h_{[n/2]}(Q)$
(follows from the Hard Lefschetz theorem for projective
algebraic varieties).
\end{enumerate}

For an arbitrary convex polytope (more generally for an Eulerian poset) 
Stanley proved ([S],thm 2.4) the property 2 above. He conjectured that 1 and 3 
also hold without the rationality hypothesis. This is still not known.

In this paper we propose an approach which we
expect to lead to a proof of 1 and 3 for general convex polytopes.
Our approach is modeled on the ``equivariant geometry'' of the (non-existent)
toric variety $X_Q$ as developed in [BL].

Namely, given a convex polytope $Q\subset\bbR^n$ we consider  the
corresponding complete fan $\Phi =\Phi _Q$ in $\bbR^n$ and
work with $\Phi$ instead of $Q$.
Let $A$ denote the graded ring of polynomial functions on $\bbR^{n}$.
Viewing $\Phi$ as a partially ordered set (of cones)
we consider a category of sheaves of $A$-modules on $\Phi$.
In this category we define a {\it minimal} sheaf $\cL_\Phi$ which corresponds
to the
$T$-equivariant intersection cohomology complex on $X_Q$ if the latter 
exists. Our first main result is the ``elementary'' decomposition theorem for 
the direct image of the minimal sheaf under subdivision of fans 
(Theorem 5.5). (Recall, that a subdivision of a fan corresponds to a proper
morphism of toric varieties.) We also develop the Borel-Moore-Verdier 
duality in the derived category of sheaves of $A$-modules on $\Phi $. 
We show that  $\cL_\Phi$ is isomorphic to its Verdier dual
(Corollary 6.26).

\remark In fact the usual (equivariant) decomposition theorem for a proper
morphism of toric varieties can be deduced from this ``elementary'' one
by the equivalence of categories proved in [L] (Thm 2.6). However the proof 
of this last result by itself uses the fundamental properties of the 
intersection cohomology. 
\endremark

The minimal sheaf $\cL_\Phi$ gives rise in a natural way to the
graded vector space $IH(\Phi)$ which we declare to be
the {\it intersection cohomology} of $\Phi $. (For rational $Q$ it is proved
in \cite{BL} that there is an isomorphism $IH(\Phi_Q)\isomo IH(X_Q)$.)
Let $ih_i(\Phi)=\dim IH^i(\Phi)$. We establish the following properties
of $IH(\Phi)$:

\begin{enumerate}
\item $\dim IH(\Phi)<\infty$;

\item $ih_i(\Phi )=0$, unless $i$ is even and $0\leq i\leq 2n$;

\item $ih_0(\Phi)=ih_{2n}(\Phi)= 1$;

\item $ih_{n-i}(\Phi)=ih_{n+i}(\Phi)$.
\end{enumerate}

Moreover, there is a natural operator $l$ of degree 2 on the space $IH(\Phi)$, 
which we expect to have the Lefschetz property as conjectured below:

\begin{conj}
For each $i\geq 1$ the map
\[
l^i:IH^{n-i}(\Phi) @>>> IH^{n+i}(\Phi)
\]
is an isomorphism.
\end{conj}

So far we were unable to prove this conjecture, but it seems to be within reach.
In case $Q$ is rational the conjecture follows from the results in [BL].
This conjecture has the following standard corollary:

\begin{cor}[of the conjecture]
$ih_i(Q)\leq ih_{i+2}(Q)$ for $0\leq i < n$.
\end{cor}

In fact the above conjecture implies ``everything'':
\begin{cor}
Assume the above conjecture is true. Then
\begin{enumerate}
\item $IH(\Phi_Q)$ is a combinatorial invariant of $Q$, i.e. it depends only
on the face lattice of $Q$. 
\item Moreover, $ih_{2j}(\Phi_Q)=h_j(Q)$ hence the $h$-vector $h(Q)$ has the 
properties conjectured by R.~Stanley.
\end{enumerate}
\end{cor}

The paper is organized as follows. 
\begin{itemize}
\item Section \ref{section:news} gives a brief account of our methods
and main results.

\item Section \ref{section:sh-hit-fan}
discusses the elementary properties of the category of abelian
sheaves on a fan and their
cohomology. 

\item In Section \ref{section:rngd-spaces} we endow a fan with the structure
of a ringed space, single out a category of sheaves of modules over the structure
sheaf and obtain the first ``geometric'' result (Theorem \ref{thm:coh-compl}).

\item In Section \ref{section:ss} we prove that our category of sheaves is
semi-simple and identify the simple objects (Theorem \ref{thm:ss}). We show that
our categories of sheaves are stable by direct image under morphisms induced by
subdivision of fans (Theorem \ref{thm:dir-im}).

\item Section \ref{section:duality} contains an account of duality on our category of sheaves. 
As a consequence we obtain the Poincare duality for $IH(\Phi )$.

\item In Section 7 we make precise Lefschetz
type conjectures and discuss their consequences.

\item In Section \ref{section:Kalai} we apply the machinery to a conjecture of
G.~Kalai (proven recently in the rational case by T.~Braden and R.D.~MacPherson)
and give our version of the proof.
\end{itemize}

\section{Summary of methods and results}\label{section:news}
\subsection{Fans as ringed spaces}
Our point of departure is the observation that a fan $\Phi$
in a (real) vector space $V$ gives rise to a topological space, which we
will denote by $\Phi$ as well, and a sheaf of graded rings $\cA_\Phi$ on it.
Namely, the points of $\Phi$ are cones, open subsets are subfans, and the
stalk $\cA_{\Phi,\sigma}$ of $\cA_\Phi$ at the cone $\sigma\in\Phi$ is
the graded algebra of polynomial functions on $\sigma$ (equivalently on the
linear span of $\sigma$) and the structure maps
are given by restriction of functions. All of these rings are quotients
of the graded algebra $A=A_V$ of polynomial functions on $V$. The grading is
assigned so that the linear functions have degree two.

In case $V$ is the Lie algebra of a torus $T$ the graded ring $A$ is
canonically isomorphic to $H^*(BT)$ -- the cohomology ring of the
classifying space of $T$.

In the case of a rational fan $\Phi $  one has the (unique) 
normal $T$-toric variety
$X_\Phi$ such that
the $T$-orbits in $X_\Phi$ are in bijective correspondence
with the cones of $\Phi$ and $\cA_{\Phi,\sigma}$ is canonically isomorphic
to the cohomology ring of the classifying space of the stabilizer of the
corresponding orbit.

All $\cA_\Phi$-modules will be regarded by default as $A$-modules. Let
$A^+$ denote the ideal of functions which vanish at the origin. For a graded
$A$-module $M$ we will denote by $\overline M$ the graded vector space
$M/A^+M$.

\subsection{A category of $\cA _{\Phi }$-modules}
To each fan $\Phi$ viewed as the ringed space $(\Phi,\cA_\Phi)$ we associate
the additive category ${\frak M}(\cA_\Phi)$ of (sheaves of finitely generated,
graded) $\cA_\Phi$-modules which are flabby and locally free over $\cA_\Phi$.
This latter condition means that, for an object $\cM$ of ${\frak M}(\cA_\Phi)$,
the stalk $\cM_\sigma$ is a free graded module of finite rank over
$\cA_{\Phi,\sigma}$. The flabbyness condition may be restated
as follows: for every cone $\sigma$ the restriction map
$\cM_\sigma @>>> \cM(\partial\sigma)$ is surjective
(where $\cM(\partial\sigma)$ is
the space of section of $\cM$ over the subfan $\partial\sigma$ consisting
of cones properly contained in $\sigma$). It is easy to see that
the sheaf $\cA_\Phi$ is flabby if and only if the fan $\Phi$ is simplicial.

In the rational case the category ${\frak M}(\cA_\Phi)$ is equivalent
to the category of equivariant perverse (maybe shifted) sheaves on $X_\Phi$.
The following theorems verify that the category ${\frak M}(\cA_\Phi)$ and the
cohomology of an object $\cM$ of ${\frak M}(\cA_\Phi)$ have the expected
properties.

Since, by definition, the objects of ${\frak M}(\cA_\Phi)$ are flabby sheaves,
it follows that, for $\cM$ in ${\frak M}(\cA_\Phi)$,
$H^i(\Phi;\cM)=0$ for $i\neq 0$.

\begin{thm}\label{thm-intro:free-coh}
Suppose that $\Phi$ is complete (i.e. the union of the cones of $\Phi$ is
all of $V$), and $\cM$ is in ${\frak M}(\cA_\Phi)$. Then, $H^0(\Phi;\cM)$
is a free $A$-module.
\end{thm}

In the rational case, $H^0(\Phi;\cM)$ is the equivariant cohomology of 
the corresponding perverse sheaf on $X_\Phi$.

The proof of Theorem
\ref{thm-intro:free-coh} rests on the observation that the cohomology
of a sheaf $\cF$ on a complete fan may be calculated by a ``cellular'' complex
$C^\bullet(\cF)$ whose component in degree $i$ is the direct sum of the
stalks of $\cF$ at cones of codimension $i$ and the differential is given by
the sum (with suitable signs) of the restriction maps. In particular, if the
sheaf $\cF$ is flabby, then the complex $C^\bullet(\cF)$ is acyclic except in
degree zero. This proves the conjecture of J.~Bernstein and the second author
(Conjecture 15.9 of \cite{BL})
on the acyclicty properties of the ``minimal complex'', which happens to
be the ``cellular complex'' of the simple object $\cL_\Phi$ (see below) of
${\frak M}(\cA_\Phi)$. In the simplicial case (when $\cL_\Phi\isomo\cA_\Phi$)
an ``elementary'' proof of this fact was given by M.~Brion in \cite{B}.

Concerning the structure of the category ${\frak M}(\cA_\Phi)$ we have the
following result.

\begin{thm}\label{thm-intro:ss}
Every object in ${\frak M}(\cA_\Phi)$ is a finite direct sum of 
indecomposable ones. The indecomposable objects are, up to a shift of
the grading, in bijective correspondence with the set of cones. 
(see Theorem 5.3 below). 
\end{thm}

\subsection{$IH$ and $IP$}
The indecomposable object of ${\frak M}(\cA_\Phi)$ which corresponds 
to the cone
$\sigma$ is a sheaf supported on the
star of $\sigma$ (which constitutes the closure of the set $\{\sigma\}$
in our topology).
Let $\cL_\Phi$ denote the indecomposable object of ${\frak M}(\cA_\Phi)$ 
which is
supported on all of $\Phi$ (the star of the origin of $V$) and whose stalk
at the origin is the one dimensional vector space in degree zero.
The fan $\Phi$ is simplicial if and only if $\cL_\Phi\isomo\cA_\Phi$.

In the rational case, when $\Phi$ is complete (and so is $X_\Phi$), the
$A$-module $H^0(\Phi;\cL_\Phi)$ is isomorphic to the $T$-equivariant
intersection cohomology $IH_T(X_\Phi)$ of $X_\Phi$ and
$\overline{IH_T(X_\Phi)}$ is the usual (non-equivariant)
intersection cohomology of $X_{\Phi }$.
This motivates the following notation:

\begin{defn}
Let $\Phi$ be a fan in $V$. Put
\[
IH(\Phi)\overset{def}{=}\overline{H^0(\Phi;\cL_\Phi)}
\]
and denote by $ih(\Phi)$ the corresponding Poincar\'e polynomial.
\end{defn}

For each cone $\sigma \in \Phi$ we may consider the corresponding
local Poincar\'e polynomial. Namely, in the rational case the graded
vector space $\overline{{\cal L}_{\Phi ,\sigma }}$  is the (cohomology
of the) stalk on the corresponding $T$-orbit $O_{\sigma }$ of the intersection
cohomology complex of $X_{\Phi }$. A normal slice to $O_{\sigma }$  is an
affine cone over some projective variety $Y_{\sigma }$. Then
$\overline{{\cal L}_{\Phi ,\sigma }}$
is the primitive part of the intersection cohomology of $Y_{\sigma }$.
This motivates the following notation.

\begin{defn}
For $\sigma \in \Phi $ put
$$IP(\sigma ):=\overline{{\cal L}_{\Phi ,\sigma }}$$
and denote by $ip(\sigma )$ the corresponding Poincar\'e polynomial.
\end{defn}

As is well known, the projectivity of a toric variety translates into
the following picture. Suppose that $\Phi$ is a complete fan in $V$ and
$l\in\cA_\Phi(\Phi)$ is a (continuous) cone-wise linear (with respect to
$\Phi$) strictly
convex function on $V$.
Multiplication by $l$ is an endomorphism (of degree 2) of ${\cal L}_{\Phi }$,
$H^0(\Phi ; {\cal L}_{\Phi })$ and $IH(\Phi )$. In the rational case it is
 the Lefschetz operator on $IH(\Phi )=IH(X_{\Phi })$ for the corresponding
projective embedding of $X_{\Phi }$. Thus we make the following
conjecture.

\begin{conj}\label{conj:lef} (Hard Lefschetz)
Let $\Phi $ be a complete fan.
Multiplication by $l$ is a Lefshetz operator on $IH(\Phi)$ i.e.
for each $i\geq 1$ the map
$$l^i:IH^{n-i}(\Phi )\to IH^{n+i}(\Phi )$$
is an isomorphism.

\end{conj}

\subsection{Subdivision and the decomposition theorem}
A fan $\Psi$ is a subdivision of a fan $\Phi$ if every cone of the latter
is a union of cones of the former. In this case there is a morphism of
ringed spaces $\pi : (\Psi,\cA_\Psi) @>>> (\Phi,\cA_\Phi)$. In the rational
case subdivision corresponds to a proper birational morphism of $T$-toric
varieties.

\begin{thm}\label{thm-intro:dir-im} (Decomposition Theorem)
The functor of direct image under subdivision restricts to the functor
$\pi_* : {\frak M}(\cA_\Psi) @>>> {\frak M}(\cA_\Phi)$.
\end{thm}

It should be pointed out that the only non-trivial part of Theorem
\ref{thm-intro:dir-im} is the fact that the direct image of a locally
free flabby sheaf is locally free which is proven by essentially the
same argument as the one used in the proof of Theorem \ref{thm-intro:free-coh}.

Combining Theorem \ref{thm-intro:dir-im} with Theorem \ref{thm-intro:ss}
we obtain the statement which  in the rational case amounts
to the Decompostion Theorem of A.~Beilinson J.~Bernstein, P.~Deligne, and
O.~Gabber (\cite{BBD}) and its equivariant analog 
(\cite{BL}) for proper birational
morphisms of toric varieties: ``the direct image of a pure object is
a direct sum of (suitably shifted) pure objects''. Continuing with notations
introduced above we have the following ``estimate'':

\begin{cor}
$\pi_*\cL_\Psi$ contains $\cL_\Phi$ as a direct summand, therefore
 $IH(\Psi)$ contains
$IH(\Phi)$ as a direct summand. So there is an inequality
\[
ih(\Psi) \geq ih(\Phi)
\]
(coefficient by coefficient) of polynomials with non-negative coefficients.
\end{cor}

\subsection{Duality}
As is well known, the (middle perversity) intersection cohomology of a
compact space admits an intersection pairing (and the same is the case
in the equivariant setting). To this end we have the following version of
Borel-Moore-Verdier duality which we develop
for the derived category of sheaves of $A$-modules on $\Phi $. 
One of the results is the following 

\begin{thm}\label{thm-intro:duality}
Let $\Phi $ be a fan in $V$. 
There is a contravariant involution $\bbD$ on ${\frak M}(\cA_\Phi)$
(i.e. a functor $\bbD : {\frak M}(\cA_\Phi)^{op} @>>> {\frak M}(\cA_\Phi)$
and an isomorphism of functors $\bbD\circ\bbD\isomo\Id$). If 
$\Phi $ is complete then there is a natural
$A$-linear non-degenerate pairing
\[
H^0(\Phi;\cM)\otimes_A H^0(\Phi;\bbD(\cM)) @>>> \omega_{A/\bbR}
\]
for every object $\cM$ of ${\frak M}(\cA_\Phi)$.
\end{thm}

Here $\omega_{A/\bbR}=A\otimes\det V^\ast$ is the dualizing $A$-module, free
of rank one, generated in degree $2\dim_\bbR V$ in accordance with our grading
convention.

If follows from Theorem \ref{thm-intro:duality} that $\bbD$ is an
   anti-equivalence
of categories, so that the dual of an indecomposable object is an 
indecomposable one.
One checks immediately that the dual $\bbD(\cL_\Phi)$ of $\cL_\Phi$
has the properties which characterize the latter. Therefore there is a
(non-canonical) isomorphism $\bbD(\cL_\Phi)\isomo\cL_\Phi$. The numerical
consequence of the auto-duality of $\cL_\Phi$ is given below.

\begin{cor}
For a complete fan $\Phi$ in a vector space of dimension $n$ the polynomial
$ih(\Phi)$ satisfies $ih_{n-k}(\Phi) = ih_{n+k}(\Phi)$. 
\end{cor}

\subsection{Kalai type inequalities}
As an application of our technology we give our restatement of the
inequality conjectured by G.~Kalai, proven in the rational case by T.~Braden
and R.D.~MacPherson in \cite{BM}.
Namely, suppose that $\Phi $ is a fan in $V$ generated by a single cone
$\sigma $, i.e. $\Phi =[\sigma ]$
and $\tau \subset \sigma$. Let $\Star (\tau )$ be the collection of cones
in $\Phi $ which contain $\tau $. Consider the minimal sheaf 
$\cL _{[\sigma ]}^{\tau }$ on $[\sigma ]$ which is based at $\tau $ (see 
5.1 below) and put 
$$IP(\Star (\tau )):=\overline{\cL _{[\sigma ], \sigma }^{\tau }}.$$

\begin{thm}
There is an inequality, coefficient by coefficient, of polynomials with
non-negative coefficients
$$ip(\sigma )\geq ip(\tau)ip(\Star(\tau)).$$
\end{thm}

\section{Abelian sheaves on fans}\label{section:sh-hit-fan}

\subsection{Fans}
A {\em fan} $\Phi$ in a real vector space $V$ of dimension $\dim V = n$
is a collection of closed convex polyhedral cones with vertex at the
origin $\Or$
satisfying
\begin{itemize}
\item any two cones in $\Phi$ intersect along a common face;
\item if a cone is in $\Phi$, then so are all of it's faces.
\end{itemize}
A fan has a structure of a partially ordered set: given cones $\sigma$ and
$\tau$ in $\Phi$ we write $\tau\leq\sigma$ if $\tau$ is a face of $\sigma$.

The origin is the unique minimal cone in every fan and will be denoted
$\Or$.

Let $d(\sigma )$ denote the dimension of the cone $\sigma $. Thus 
$d(\sigma )=0$ iff $\sigma =\Or $. 

The fan $\Phi$ is {\em complete} if and only if the union of cones of $\Phi$
is all of $V$.
The fan $\Phi$ is {\em simplicial} if every cone
of $\Phi$ is simplicial. A cone of dimension $k$ is {\em simplicial}  if it has
$k$ one-dimensional faces (rays).

\subsection{Topology on a fan}
The (partially ordered) $\Phi$ will be considered as a topological space
with the open sets the subfans of $\Phi$. A subset, say $S$, of $\Phi$
generates a subfan, denoted $\lbrack S\rbrack$.

An open subset is {\em irreducible} if it is not a union of two open subsets
properly contained in it. The irreducible open sets are the subfans generated
by single cones. We will frequently abuse notation and write $[\sigma]$ for
the irreducible open set $\lbrack\{\sigma\}\rbrack$.

Note that the topological space $\Phi$ has the following property:
the intersection of irreducible open sets is irreducible.

Let $\Phi_{\leq k}$ denote the subset of cones of dimension at most $k$;
this is an open subset of $\Phi$.

Let $\sigma $ be a cone in $\Phi $. Denote by $\Star(\sigma)\subset \Phi$
the subset of all the cones $\tau $ such that $\sigma \leq \tau$. This is a
closed subset of $\Phi $. Its image under the projection
$V\to V/\Span(\sigma )$ is a fan that will be denoted by
$\overline{\Star(\sigma )}$.

\subsection{Sheaves on a fan}
Regarding $\Phi$ as a topological space with open sets the subfans of $\Phi$,
we consider sheaves on $\Phi$.
 Let $\cI(\Phi)$ denote the partially ordered
set of irreducible open sets of $\Phi$ and inclusions
thereof. This partially ordered set is isomorphic to $\Phi$.

A sheaf on $\Phi$ restricts to a presheaf (a contravariant functor) on
$\cI(\Phi)$ and this correspondence is an equivalence of categories.
Since, for a sheaf $F$ and a cone $\sigma$, the stalk $F_\sigma$ is
equal to the sections $\Gamma([\sigma];F)$ of $F$ over the
corresponding irreducible open set, the sheaf $F$ is uniquely determined
by the assignment $\sigma\mapsto F_\sigma$ and the restriction maps
$F_\sigma @>>> F_\tau$ whenever $\tau\leq\sigma$.

As usual the support $\supp (F)$ of a sheaf $F$ is the closure of the set 
of $\sigma $'s, such that $F_{\sigma }\neq 0$.

For a cone $\sigma$ let $\partial\sigma$ denote the subfan generated by the
proper faces of $\sigma$. Let $F$ be a sheaf and consider the following
condition:
\begin{multline}\label{cond-fl}
\text{for every cone $\sigma$ the canonical map} \\
\text{$F_\sigma @>>> \Gamma(\partial\sigma;F)$ is surjective.}
\end{multline}

\begin{lemma}\label{lemma:flabby}
A sheaf $F$ satisfies the condition \eqref{cond-fl} if and only if it
is flabby.
\end{lemma}
\begin{pf}
The condition is obviously necessary. To see that it is sufficient we
need to show that a section of $F$ defined over a subfan $\Psi$ extends
to a global section. Clearly, it is sufficient to show that a section
defined over $\Psi\bigcup\Phi_{\leq k}$ extends to a section defined over
$\Psi\bigcup\Phi_{\leq k+1}$, but this is immediate from \eqref{cond-fl}.
\end{pf}

For the rest of this section 3 we restrict our
attention to the category of
sheaves of $\bbR $-vector spaces on $\Phi$ which we will denote by $\Sh(\Phi)$.

\subsection{Some elementary properties of the category $\Sh(\Phi)$}
For a cone $\sigma$ let $i_\sigma : \{\sigma\}\hookrightarrow\Phi$
denote the inclusion. The embedding $i_\sigma$ is locally closed, and
closed (respectively open) if and only if $\sigma$ is maximal (respectively
minimal, i.e. the origin).

Suppose that $W$ is a vector space considered as a sheaf on $\{\sigma\}$.
Then, clearly, $\R^p(i_\sigma)_*W = 0$ for $p\neq 0$, the sheaf
$(i_\sigma)_*W$ is an injective object in $\Sh(\Phi)$ equal to the constant
sheaf $W$ supported on $\Star(\sigma)$.

Let $i_{[\sigma]} : [\sigma]\hookrightarrow\Phi$ denote the open embedding
of the irreducible open set $[\sigma]$. Then, the extension by zero
$(i_{[\sigma]})_!W$ of the constant sheaf $W$ on $[\sigma]$ is a projective
object in $\Sh(\Phi)$. Thus the abelian category $\Sh(\Phi )$ has enough
projectives.

\subsection{The cellular complex of a sheaf}
Let $\Phi $ be a fan in $V\simeq {\Bbb R}^n$ and $F$ be a sheaf on $\Phi $.
Choose an orientation of each cone in $\Phi $. The
{\it cellular complex} $C^{\bullet}(F)$ of $F$ is defined as 
follows
$$C^{\bullet}(F):=0\to C^0(F)\to C^1(F) \to ...\to C^n(F)\to 0,$$
where
$$C^i(F)=\underset{d(\sigma)=n-i}{\bigoplus }F_{\sigma}$$
and the differential $d^i:C^i(F)\to C^{i+1}(F)$ is the sum of the
restriction maps $F_{\sigma }\to F_{\tau }$ with the sign $\pm 1$
depending on whether the orientations of $\sigma $ and $\tau $
agree or disagree.

\begin{remark} $C^{\bullet}(\cdot )$ is an exact functor from $\Sh(\Phi )$ to
complexes of vector spaces.
\end{remark}

\begin{prop}
Assume that the fan $\Phi $ is complete. Then the complex $C^{\bullet}(F)$
 is quasi-isomorphic to $\R \Gamma (\Phi ;F)$. In particular,
$$H^i(C^{\bullet}(F))=H^i(\Phi; F).$$
\end{prop}

\begin{pf}
First notice that since $\Phi $ is complete we have
$$\Gamma (\Phi ;F)=H^0(C^{\bullet}(F)).$$
Indeed, consider the open covering of $\Phi $ by maximal irreducible open
sets $[\sigma ]$, $d(\sigma )=n$. Then to give a global section $s\in
\Gamma (\Phi ;F)$ is the same as to give a collection of local
sections $s_{\sigma }\in F([\sigma ])$, $d(\sigma )=n$ such that
$s_{\sigma }=s_{\tau }$ in $F([\sigma \cap \tau ])$ if $d(\sigma \cap \tau )=
n-1$.
This shows that $\Gamma (\Phi ;F)=H^0(C^{\bullet}(F))$.

To prove the proposition consider an injective resolution
$$F\to I^0\to I^1\to ...,$$
where sheaves $I^j$ are direct sums of constant sheaves on closed subsets
$\Star(\sigma )$, $\sigma \in \Phi$ extended by zero to $\Phi $. Since
$C^{\bullet}(\cdot)$ is an exact functor it suffices to prove the following
claim.

\begin{claim}
Let $\Psi $ be a fan in $V$, $\sigma \in \Psi$ and $W$ -- a vector space.
Assume that the fan $\overline{\Star(\sigma)}$
in $V/\Span(\sigma)$ is complete.
 Let $I$ be the constant
sheaf $W$ on $\Star(\sigma )$ extended by zero to $\Psi $. Then
$$H^i(C^{\bullet}(I))=\left\{
\begin{array}{cl}
W&if\ \ i=0\\
0,&otherwise.\\
\end{array}
\right.
$$
\end{claim}

\noindent{\it Proof of claim.} Since $\Star(\sigma )$ is a closed subset
of $\Psi $ which is isomorphic to the complete fan $\overline{\Star(\sigma )}$
  (of dimension $n-d(\sigma)$) we may assume that $\sigma =\Or$ and hence
$\Star(\sigma )=\Psi $ is a complete fan in $V$. 
Then $C^{\bullet}(I)$ is isomorphic to a cellular cochain complex of a ball 
of dimension $n$ (with coefficients $W$).
This proves the claim and the proposition.

\end{pf}

\begin{cor}
Let $\sigma \subset V$ be a cone and $\Psi =[\partial \sigma ]$ be the fan
generated by the boundary $\partial \sigma $. Then for any sheaf $F$ on $\Psi $
 the shifted cellular complex $C^\bullet (F)[n-d(\sigma )+1]$ is quasi-isomorphic
to $\R \Gamma (\Psi ;F)$.
\end{cor}

\begin{pf}
This follows from the above proposition since the fan $\Psi $ is isomorphic
to a complete fan in ${\Bbb R}^{d(\sigma )-1}$.
\end{pf}

 Later on we will need the following version of the previous proposition.

\begin{prop}
Let $\Phi $ be a fan in $V$, not necessarily complete. Let $F\in \Sh(\Phi)$
be such that its support $Z=Supp(F)$ satisfies the following condition:
 for each $\sigma \in Z$ the fan $\overline{\Star(\sigma)}$
in $V/\Span(\sigma )$ is complete. Then the cellular
complex $C^{\bullet}(F)$ is quasi-isomorphic to $\R \Gamma(\Phi;F)$.
\end{prop}

\begin{pf} Same as that of Proposition 3.3.
\end{pf}

The next three lemmas will be used later on.

\begin{lemma}
Let $\sigma \subset V$ be a cone of positive dimension, i.e. 
$\sigma \neq \Or $, and $\Phi =[\sigma ]$ be the fan generated
by $\sigma $. Let $F$ be a constant sheaf on $\Phi $. Then the cellular
complex $C^\bullet (F)$ is acyclic.
\end{lemma}

\begin{pf}
We may assume that $d(\sigma )=dimV=n$. Then the cellular complex
$C^\bullet (F)$ is isomorphic to an augmented chain complex of
a ball of dimension $n-1$. This proves the lemma.
\end{pf}

\begin{remark}
In the notation of the previous lemma notice that if $F\neq 0$, then
$H^0([\sigma ];F)\neq 0$. Thus for a fan which is not complete the
cellular complex does not necessarily compute the cohomology of the
sheaf.
\end{remark}

\begin{lemma}
Let $\sigma $ be a cone in $V$ and consider the fan $[\sigma ]$. 
Let $F$ be a flabby
sheaf on the fan $[\sigma ]$. Then the cellular complex $C^\bullet (F)$ is
acyclic except in the lowest degree, i.e. $H^i(C^\bullet (F))=0$ for
$i\neq n-d(\sigma )$.
\end{lemma}

\begin{pf}
Put $\Psi :=[\partial \sigma ]$ and denote $F_{\Psi }:=F\vert _{\Psi }$. 
Then by Corollary 3.5 
$H^{i+(n-d(\sigma )+1)}(C^\bullet (F_{\Psi }))=H^i(\Psi ;F)$.
 Since $F_{\Psi }$ is flabby
$H^j(C^\bullet (F_{\Psi }))=0$ if $j\neq n-d(\sigma )+1$, and
$H^{n-d(\sigma )+1}(C^\bullet (F_{\Psi }))=\Gamma (\Psi ;F)$. Again by the
flabbiness of $F$ the map $F_{\sigma }\to \Gamma (\Psi ;F)$ is surjective,
hence $H^i(C^\bullet (F))=0$ for $i\neq n-d(\sigma )$.
\end{pf}

\subsection{Cohomology of some simple sheaves}
\begin{lemma}
Let $\sigma \subset V$ be a cone and consider the fan $\Phi =[\sigma ]$.
Then for any sheaf $F$ on $\Phi $ $H^i(\Phi ,F)=0$, if $i>0$.
\end{lemma}

\begin{pf}
Let ${\Bbb R}_{\Phi }$ be the constant sheaf on $\Phi $ with stalk ${\Bbb R}$.
We have the isomorphism of functors
$$\Gamma (\Phi , \cdot )=\Hom ({\Bbb R}_{\Phi }, \cdot ).$$
But ${\Bbb R}_{\Phi }$ is a projective object in $Sh(\Phi )$.
 This proves the lemma.
\end{pf}

\begin{lemma}
Let $\Phi $ be a fan in $V$ and $F\in Sh(\Phi )$ be the constant sheaf on
$\Phi $ with stalk $W$. Then
$$H^i(\Phi ,F)=\left\{
\begin{array}{ll}
W, & if\ \ i=0\\
0, & otherwise.\\
\end{array}
\right.
$$
\end{lemma}

\begin{pf}
Since the space $\Phi $ is connected $H^0(\Phi ,F)=W$. The rest follows from
the injectivity of the sheaf $F$.
\end{pf}

\begin{lemma}
Let $\Phi $ be a complete fan in $V$ and $W$ be a vector space. Consider
the sheaf $W_{\Or }$ on $\Phi $ which is the extension by zero of the
sheaf $W$ on the open point $\Or $. Then
$$H^i(\Phi ,W_{\Or })=\left\{
\begin{array}{ll}
W, & if\ \ i=n\\
0, & otherwise.\\
\end{array}
\right.$$
\end{lemma}

\begin{pf}
Put $Z:=\Phi - \{\Or \}$ and let $i:Z\hookrightarrow \Phi $ be the corresponding
closed embedding. Let $W_{\Phi }$, $W_Z$ denote the constant sheaves with
 stalk $W$ on $\Phi $ and $Z$ respectively. We have the exact sequence
$$0\to W_{\Or }\to W_{\Phi }\to i_*W_Z\to 0.$$
By Lemma 3.11 $H^0(\Phi ,W_{\Phi })=W$ and $H^i(\Phi ,W_{\Phi })=0$ if $i>0$.
On the other hand the cellular complex $C^\bullet (i_*W_Z)$ is isomorphic
to a chain complex of a sphere of dimension $n-1$. Hence by Proposition 3.6 
$$H^j(\Phi ,i_*W_Z)=\left\{
\begin{array}{ll}
W, & if\ \ j=0,n-1\\
0, & otherwise\\
\end{array}
\right.
$$
if $n>1$ and
$$H^j(\Phi ,i_*W_Z)=\left\{
\begin{array}{ll}
W\oplus W, & if\ \ j=0\\
0, & otherwise\\
\end{array}
\right.
$$
if $n=1$. The map $H^0(\Phi ,W_{\Phi })\to H^0(\Phi ,i_*W_Z)$ is injective,
so the lemma follows.
\end{pf}

\subsection{The cocellular complex}
Let $\Phi $ be a fan in $V$ and $\bbR _{\Or }$ be the extension by zero
to $\Phi $ of the sheaf $\bbR$ on the open point $\Or \in \Phi$.
We want to construct an injective resolution of $\bbR _{\Or }$ of a special
kind. For a cone $\sigma \in \Phi $ let $i_{\sigma }:\{ \sigma \}
\hookrightarrow \Phi $ denote the inclusion. Consider the sheaf
$i_{\sigma *}\bbR $ on $\Phi $. This is a constant sheaf on $\Star (\sigma )$
extended by zero to $\Phi $.
Such sheaves are injective. Note that for $\tau \leq \sigma $ we have
the natural surjective morphism $r_{\tau \sigma }:i_{\sigma *}\bbR \to
 i_{\tau *}\bbR$. Consider the complex
$$K^\bullet =K^\bullet _{\Phi }:=K^{-n}\overset{d^{-n}}{\rightarrow }K^{-n+1}
\overset{d^{-n+1}}{\rightarrow }...\overset{d^{-1}}{\rightarrow}K^0,$$
where
$$K^{-n+j}=\underset{d(\sigma )=j}{\bigoplus }i_{\sigma *}\bbR,$$
and the differential $d$ is the sum of maps $r_{\tau \sigma }$ with
$\pm $ sign depending on the compatibility of orientations of $\sigma $ and
$\tau $.

Note that $K^{-n}$ is the constant sheaf $\bbR$ on $\Phi $ and $\Ker d^{-n}=
\bbR _{\Or }$.

\begin{lemma}
The complex $K^\bullet $ is an injective resolution of the sheaf
 $\bbR _{\Or }[n]$.
\end{lemma}

\begin{pf}
The injectivity of $K^\bullet$ was noted before. It remains to show that
for $\Or \neq \sigma \in \Phi $ the complex of stalks
$$K^{-n}_{\sigma }\to K_{\sigma }^{-n+1} \to ... \to K_{\sigma }^0$$
is exact. Since $i_{\tau *}\bbR$ is the constant sheaf on $\Star (\tau )$
the above complex is isomorphic to the  dual of the cellular
complex $C^\bullet (\bbR _{[\sigma ]})$.
This last complex is acyclic by Lemma 3.7 above.

\end{pf}

\begin{defn}
Given a sheaf $F$ of $\bbR $-vector spaces on $\Phi $ we define its
{\it cocellular complex} $K^\bullet (F)=K^\bullet _{\Phi }(F)$ as
$$K^\bullet (F):=\Hom ^\bullet (F,K^\bullet ).$$
\end{defn}

\begin{remark} The cocellular complex has the following properties.
\begin{enumerate}
\item $K^i (F)=0$ unless $i\in [-n,0]$.

\item By the above lemma $K^\bullet (F)=\R \Hom ^\bullet (F,\bbR _{\Or }[n])$.

\item Note that $\Hom (F, i_{\sigma *}\bbR )=\Hom (F_{\sigma },\bbR )$.
Thus the cocellular complex $K^\bullet (F)$ is just the
 dual of the cellular complex $C^\bullet (F)$:
$$K^\bullet (F)=\Hom ^\bullet (C^\bullet (F),\bbR ).$$

\item An appropriate version of the cocellular complex will play a role
in our discussion of the Borel-Moore-Verdier duality for sheaves of
$A_{\Phi }$-modules (see Section 6) on complete fans. Of course, similar
duality exists for sheaves of vector spaces, but we will not need it.
\end{enumerate}
\end{remark}

\section{Fans as ringed spaces}\label{section:rngd-spaces}
\subsection{The structure sheaf of a fan}
Let $\Phi $ be a fan in $V$. 
Let $V^\ast$ denote the constant sheaf $\sigma\mapsto V^\ast$ on $\Phi$.
Let $\Omega^1_\Phi$ denote the subsheaf of $V^\ast$ given by
\[
\Phi\ni\sigma\mapsto\Omega^1_{\Phi,\sigma}
\overset{def}{=}\sigma^\perp\subseteq V^\ast\ ,
\]
where $\sigma^\perp$ denotes the subspace of linear functions which vanish
identically on $\sigma$.

Let $\cG$ denote the sheaf determined by the assignment
\[
\Phi\ni\sigma\mapsto\cG_\sigma = \Span(\sigma)^\ast\ .
\]
Thus, there is a short exact sequence of sheaves
\begin{equation}\label{ses:Omega}
0 @>>> \Omega^1_\Phi @>>> V^\ast @>{\pi}>> \cG @>>> 0\ .
\end{equation}

>From now on $A$ will denote the symmetric algebra of $V^\ast$ with grading
determined by assigning degree $2$ to $V^\ast$. We will use the notation
$A_{\Phi }$ for the corresponding constant sheaf on $\Phi $.

\begin{defn}
The structure sheaf $\cA_\Phi$ is the symmetric algebra of $\cG$, i.e.
the sheaf of cone-wise polynomial functions, graded so that the
linear functions have degree 2.
\end{defn}

\begin{remark}
Clearly, there is an epimorphism of sheaves of
graded algebras $A_{\Phi } @>>> \cA_\Phi$.
\end{remark}
\begin{remark}
With these definitions $(\Phi,\cA_\Phi)$ is a ringed space over the
one point ringed space $(\emptyset,A)$ which we imagine as ``the empty fan'' in
$V$.
\end{remark}

In what follows ``an $\cA_\Phi$-module''
will mean ``a (locally) finitely generated graded $\cA_\Phi$-module''
 and similarly
for $A_{\Phi }$-modules. An $\cA_\Phi$-module $\cM$ is {\it locally free} 
if, for every
cone $\sigma$, $\cM_\sigma$ is a free (graded) $\cA_{\Phi,\sigma}$-module.

Let $A^+$ denote the ideal of elements of positive degree. For an $A$-module
$M$ we will denote by $\overline{M}$ the graded vector space
$M/MA^+$.

For a graded $A$-module (or sheaf) $M=\oplus M_k$ denote by
$M(t)$ the corresponding shifted object
$$M(t)_k=M_{k+t}.$$

The flabbiness criterion (3.1) applied to an $\cA_\Phi$-module together with
Nakayama's Lemma amounts to the following.

\begin{lemma}
An $\cA_\Phi$-module $\cM$ is flabby if and only if for every cone $\sigma$
the canonical map
\[
\overline{\cM_\sigma} @>>> \overline{\Gamma(\partial\sigma;\cM)}
\]
is surjective.
\end{lemma}

\begin{defn}
Let ${\frak M}(\cA_\Phi)$ denote the additive category of flabby locally free
$\cA_\Phi$-modules considered as a full subcategory of sheaves of
$\cA_\Phi$-modules and morphisms of dergee zero.
\end{defn}

\begin{lemma}
The fan $\Phi$ is {\em simplicial} if and only if
the structure sheaf $\cA_\Phi$ is flabby.
\end{lemma}
\begin{pf}
The statement is easily seen to be equivalent to the following one:
{\em Suppose that $P_1,\dots,P_n$ are polynomials in variables $x_1,\dots,
x_n$  \\ which satisfy $\displaystyle\frac{\partial P_i}{\partial x_i} = 0$ and
$P_i\vert_{x_j=0} = P_j\vert_{x_i=0}$ for all $i$ and $j$. Then there is a
polynomial $Q$ such that $Q\vert_{x_i=0} = P_i$.}
Verification of the latter fact is left to the reader.
\end{pf}

\subsection{Cohomology of objects in ${\frak M}(\cA_\Phi )$ on complete fans}

\begin{thm}\label{thm:coh-compl}
Suppose that $\Phi$ is a complete fan in $V$ and $\cM$ is in
${\frak M}(\cA_\Phi)$. Then, $H^i(\Phi;\cM) = 0$ for $i\neq 0$ and
$H^0(\Phi;\cM)$ is a free $A$-module.
\end{thm}
\begin{pf}
The vanishing of higher cohomology is a direct consequence of the flabbiness
of the objects of ${\frak M}(\cA_\Phi)$.

A graded $A$-module
$M$ is free if and only if $\Ext^i_{A}(M,A)=0$
for $i\neq 0$.
By Proposition 3.3 above $H^0(\Phi ;\cM )=H^0C^\bullet (\cM )$.
Note that, for every $j$,
$\Ext^i_{A}(C^j(\cM),A)=0$ for $i\neq j$. Since, in addition,
$H^iC^\bullet(\cM) = 0$ for $i\neq 0$ it follows by the standard argument
that $\Ext^i_{A}(H^0C^\bullet(\cM),A)=0$ for $i\neq 0$.
\end{pf}

\section{Minimal sheaves and intersection cohomology}\label{section:ss}
\subsection{Minimal sheaves}
Recall that, for a cone $\sigma$, the set $\Star(\sigma)$ is defined as the
collection of those cones $\tau$ which satisfy $\tau\geq\sigma$. Namely,
$\Star(\sigma)$ is the closure of the set $\{\sigma\}$.

Fix a cone $\sigma$ and consider the following conditions on an object $\cM$
of ${\frak M}(\cA_\Phi)$:
\begin{enumerate}
\item $\cM_{\sigma }\neq 0$ and $\cM_\tau\neq 0$ only if
$\tau\in\Star(\sigma)$;

\item for every cone $\tau\in\Star(\sigma)$, $\tau\neq\sigma$, the canonical
map
$\overline{\cM_\tau} @>>> \overline{\Gamma(\partial\tau;\cM)}$
is an isomorphism.
\end{enumerate}

\begin{defn}
In what follows we will refer to an object as above as {\em a minimal
sheaf based at $\sigma$}. An instance of a minimal sheaf $\cM$ based at
$\sigma$ with $\cM_\sigma\isomo\cA_{\Phi,\sigma}$ will be denoted
$\cL_\Phi^\sigma$. We will also denote $\cL_\Phi=\cL_\Phi^{\Or}$
\end{defn}

\begin{prop}
Let $\Phi $ be any fan.

\begin{enumerate}
\item For every $\sigma \in \Phi$ and every finitely generated graded
$\cA_{\Phi, \sigma}$-module $M$ there exists a unique (up to an isomorphism)
minimal sheaf $\cM $ based at $\sigma $ such that $\cM_{\sigma }=M$.
In particular, the minimal sheaf $\cL_{\Phi }^{\sigma}$ exists for each
$\sigma \in \Phi$.
\item Moreover, if $\cM $ is a minimal sheaf based at $\sigma $
then
$$\cM \simeq \cL _{\Phi }^{\sigma }\otimes _{\Bbb R}\overline{\cM_{\sigma }}.$$
In particular $\cM $ is a direct sum of sheaves $\cL _{\Phi }^{\sigma }(t)$, 
$t\in \bbZ $. 
\item The minimal sheaves $\cL _{\Phi }^{\sigma }(t)$, $t\in \bbZ$ are 
indecomposable objects in the category of $\cA _{\Phi }$-modules.
\item Let $U\subset \Phi $ be an open subset, i.e. $U$ is a subfan. Then 
$$\cL _{\Phi }\vert _U=\cL _U.$$

\end{enumerate}
\end{prop}

\begin{pf}
Easy exercise.
\end{pf}

\begin{thm}\label{thm:ss}
Let $\Phi $ be a fan in $V$. 
Every  object  $\cM $ of
${\frak M}(\cA_\Phi)$ is
 isomorphic to a
direct sum of minimal sheaves. In particular, $\cM $ is a direct sum 
of indecomposable objects $\cL _{\Phi }^{\sigma }(t)$, $t\in \bbZ$.
\end{thm}
\begin{pf}
The last part of the theorem follows from the first one using parts 2 and 3 
of Proposition 5.2. 

Consider an object $\cM$ of ${\frak M}(\cA_\Phi)$. We will show that it is
isomorphic to a direct sum of minimal sheaves by induction on
$$\vert\cM\vert=\displaystyle\sum_{\sigma}rank_{\cA_{\Phi,\sigma}}\cM_\sigma .$$

Consider a cone $\sigma$ such that $\cM_\sigma\neq 0$ and,
for every $\tau\in\partial\sigma$, $\cM_\tau = 0$. We will show that,
for each such $\sigma$, $\cM$ contains as a direct summand a minimal sheaf
$\cK$ based at $\sigma$ with $\cK_\sigma = \cM_\sigma$. That is, we will
construct a direct sum decompostion
\begin{equation}\label{dir-sum}
\cM = \cK\oplus\cN
\end{equation}
with $\cK$ as above and $\cN$ in ${\frak M}(\cA_\Phi)$. This is sufficient,
since, clearly, $\vert\cN\vert < \vert\cM\vert$.

Note, that we need to specify the direct sum decompostion \eqref{dir-sum} only
on $\Star(\sigma)$. Therefore it is sufficient to treat the case when
$\sigma$ is the origin $\Or$.

We proceed to construct the decomposition \eqref{dir-sum} by
induction on the dimension of the cone and the number of cones of the
given dimension.

Let $\cK_\Or = \cM_\Or$ and $\cN_\Or = 0$.
Assume that a direct sum decomposition
\[
\cM\vert_{\Phi_{\leq k}} = \cK_{\leq k}\oplus\cN_{\leq k}
\]
in ${\frak M}(\cA_\Phi\vert_{\Phi_{\leq k}})$ has been defined and consider
a cone $\sigma$ of dimension $k+1$. Since $\partial\sigma$ consists
of cones of dimension at most $k$ the induction hypothesis
says that there is a direct sum decompostion of $\Gamma (\partial \sigma ; 
\cA _{\Phi })$-modules
\[
\Gamma(\partial\sigma;\cM) =
\Gamma(\partial\sigma;\cK_{\leq k})\oplus
\Gamma(\partial\sigma;\cN_{\leq k})\ .
\]
\begin{claim}
There is a decomposition
\[
\cM_\sigma = K_\sigma\oplus N_\sigma
\]
into a direct sum of free $\cA_{\Phi,\sigma}$-modules, such that
the restriction homomorphism
$$\cM_\sigma \twoheadrightarrow \Gamma(\partial\sigma;\cM)$$
maps $K_{\sigma }$ to $\Gamma(\partial\sigma;\cK_{\leq k})$ and
$N_{\sigma }$ to $\Gamma(\partial\sigma;\cN_{\leq k})$ and induces
an isomorphism
$$
\overline{K_\sigma}\stackrel{\sim}{\rightarrow}
\overline{\Gamma(\partial\sigma;\cK_{\leq k})} $$
and an epimorphism
$$
\overline{N_\sigma}
\twoheadrightarrow
\overline{\Gamma(\partial\sigma;\cN_{\leq k})}$$
\end{claim}

Assume the claim for the moment. The desired extension $\cK_{\leq k+1}$
(respectively $\cN_{\leq k+1}$) of $\cK_{\leq k}$ (respectively $\cN_{\leq k}$)
is given, for every cone $\sigma$ of dimension
$k+1$, by $\cK_{\leq k+1,\sigma} = K_\sigma$ (respectively
$\cN_{\leq k+1,\sigma} = N_\sigma$) and has all the required properties.

\noindent{\it Proof of claim.}
Choose a subspace $Z\subset \Gamma(\partial\sigma;\cK_{\leq k})$ which maps
isomorphically onto $\overline{\Gamma(\partial\sigma;\cK_{\leq k})}$ under the
residue map. Choose a subspace $S\subset\cM_\sigma$ so that the map
$\cM_\sigma @>>> \Gamma(\partial\sigma;\cM)$ restricts to an isomorphisms
$S \stackrel{\sim}{\rightarrow} Z$. Since $S\bigcap A^+\cM_\sigma = 0$ there is a subspace
$T\subset\cM_\sigma $ such that $S\bigcap T=0$ and $S\oplus T$ generates
$\cM_\sigma $ freely. Subtracting, if necessary, elements of $A\cdot S$ from
 elements
of $T$ we may assume that the image of $T$ under the map
$\cM_\sigma @>>> \Gamma(\partial\sigma;\cM)$ is contained in
$\Gamma(\partial\sigma;\cN_{\leq k})$. Thus we may take
$K_\sigma = \cA_{\Phi,\sigma}S$ and $N_\sigma = \cA_{\Phi,\sigma}T$.
This concludes the proof of the claim and of the theorem.
\end{pf}

\subsection{Subdivision of fans and the decomposition theorem}
Suppose that $\Phi$ and $\Psi$ are two fans in $V$ and $\Psi$ is a subdivision
of $\Phi$, which is to say, every cone of $\Phi$ is a union of cones of
$\Psi$. (In the rational case this induces a proper morphism of toric varieties).
 This corresponds to a morphism of ringed spaces
$\pi : (\Psi,\cA_\Psi) @>>> (\Phi,\cA_\Phi)$.
 The next theorem combined
with the structure Theorem 5.3 is a combinatorial
 analog of
the decomposition theorem ([BBD],[BL]).

\begin{thm}\label{thm:dir-im}
In the notations introduced above, for $\cM$ in ${\frak M}(\cA_\Psi)$,
\begin{enumerate}
\item $\R^i\pi_*\cM = 0$ for $i\neq 0$ and $\pi_*\cM$ is flabby;

\item $\pi_*\cM$ is locally free.
\end{enumerate}
In other words, the direct image under subdivision $\pi$ restricts to an exact
functor $\pi_* : {\frak M}(\cA_\Psi) @>>> {\frak M}(\cA_\Phi)$.
\end{thm}
\begin{pf}
The first claim follows from the flabbiness of $\cM$.

Since the issue is local on $\Phi$ we may assume that the latter is generated
by a single cone $\sigma$ of top dimension $n$, i.e. $\Phi =[\sigma ]$.
By induction on dimension it is
sufficient to show that the stalk
$$(\pi_*\cM)_\sigma = H^0(\Phi;\pi_*\cM) = H^0(\Psi ;\cM )$$ 
is a free $A$-module.

Let $Z=\pi^{-1}(\sigma)$. This is a closed subset of $\Psi$ which consists
of the cones which subdivide the interior of $\sigma$.

\begin{claim} For any sheaf $F$ on $\Psi $ the restriction map 
$H^0(\Psi ;F)\to H^0(Z;F)$ is an isomorphism.
\end{claim}

\medskip
\noindent{\it Proof of claim.} Indeed, a global section $\alpha \in 
\Gamma (\Psi ;F)$ is the same as a collection of local sections 
$\alpha _{\tau }\in F_{\tau }=\Gamma ([\tau ];F)$, $d(\tau )=n$ such that 
$\alpha _{\tau }=\alpha _{\xi }$ in $F_{\tau \cap \xi }$ in case 
$d(\tau \cap \xi )=n-1$. The same local data specifies an element 
in $\Gamma (Z;F)$. This proves the claim.
\medskip

Denote by $i:Z\hookrightarrow \Psi$ the corresponding closed embedding.
Put
$$\cM_Z:=i_*i^*\cM.$$
Then by the above claim $(\pi _*\cM)_{\sigma }=H^0(\Psi; \cM_Z).$

By Proposition 3.6 the cellular complex $C^{\bullet}(\cM_Z)$ is
quasi-isomorphic to $\R \Gamma(\Psi;\cM_Z)$. Note that the sheaf $\cM_Z$ is
flabby. Now we show that the $A$-module $H^0(\Psi;\cM_Z)$ is free by the
same argument as in the proof of Theorem 4.7.
\end{pf}

\subsection{Intersection cohomology of fans}
\begin{defn}
Let $\Phi $ be a fan in $V$. Define its {\it intersection cohomology}
as the graded vector space
$$IH(\Phi ):=\overline{H^0(\Phi ;\cL _{\Phi })}$$
and denote by $ih(\Phi )$ the corresponding Poincare polynomial.
\end{defn}

\begin{lemma}
1. $\dim IH(\Phi )<\infty .$

2. $ih_j(\Phi )=0$ for $j<0$ or $j$ odd.

3. $ih_0(\Phi )=1$.
\end{lemma}

\begin{pf}
The $A$-module $H^0(\Phi ;\cL _{\Phi })$ is finitely generated, so 1) follows.
The definition of $\cL _{\Phi }$ implies that for each $\sigma
\in \Phi $ the (graded) $A$-module $\cL _{\Phi ,\sigma }$ has no negative or 
odd 
part
 (use induction on $d(\sigma )$). Thus the same is true for
$H^0(\Phi; \cL _{\Phi })$ and 2) follows. Also by induction on $d(\sigma )$
(and using Lemma 3.11 above) one checks that the zero component of 
$\cL _{\Phi ,
\sigma }$ has dimension 1. In other words the zero component of the sheaf
$\cL _{\Phi }$ is the constant sheaf ${\Bbb R}$. Thus by Lemma 3.11 
$ih _0(\Phi )=1$.
\end{pf}

We will be able to say much more about $IH(\Phi )$ in case $\Phi $ is a
{\it complete} fan.

\begin{lemma}
Let $\Psi $ be a subdivision of $\Phi $. Then $IH(\Phi )$ is a direct summand
of $IH(\Psi )$. In particular, one has the inequality
$$ih(\Psi )\geq ih(\Phi )$$
(coefficient by coefficient) of polynomials with nonnegative coefficients.
\end{lemma}

\begin{pf}
Let $\pi :(\Psi ,\cA _{\Psi })\to (\Phi ,\cA _{\Phi })$ denote the
corresponding morphism of ringed spaces. Then $(\pi _*\cL _{\Psi })_{\Or }=
{\Bbb R}$. Hence the sheaf $(\pi _*\cL _{\Psi })$ contains $\cL _{\Phi }$
as a direct summand (see Theorem 5.3). The lemma follows.
\end{pf}

\section{Borel-Moore-Verdier duality}\label{section:duality}
Let $\Phi $ be a fan in $V={\Bbb R}^n$.
 Let $A_{\Phi }$ as usual be
the constant sheaf on $\Phi $ with stalk $A$. Denote by
$D^b_c(A_{\Phi }-mod)$ the bounded derived category of (locally
finitely generated) $A_{\Phi }$-modules. In particular the additive category
of sheaves
${\frak M}(\cA_\Phi)$ is a full subcategory of $D^b_c(A_{\Phi}-mod)$.

 In this section we define the duality functor,
i.e. a contravariant involution $\bbD$ on the category $D^b_c(A_{\Phi}-mod)$.
 The duality preserves the subcategory
 ${\frak M}(\cA_\Phi)$ and $\bbD(\cL_\Phi)\simeq \cL_\Phi$.
Among other things, it gives rise to Poicar\'e duality in $IH(\Phi)$.

The contsruction of $\bbD$ follows the general pattern of Borel-Moore-Verdier
 duality.
In particular, for any sheaf on $\Phi$ we define explicitly the co-sheaf of
''compactly supported`` sections.

\subsection{Co-sheaves and homology}
While a sheaf on a fan $\Phi$ (with values in an abelian category $\cC$)
 is a functor
$\Phi^0 @>>> \cC$, a {\em co-sheaf} is a functor $\Phi @>>> \cC$. A 
co-sheaf determines 
a functor $\Phi^0 @>>> \cC^0$. 
Given a sheaf $F$ on $\Phi $ its space of global sections $\Gamma (\Phi; F)$
is by definition the inverse limit
$\underset{\underset{\Phi}{\leftarrow}}{\lim }F$.
This is a left exact functor. For a co-sheaf $\cV: \Phi \to \cC$ we define its
space of global co-sections as the direct limit
$\underset{\underset{\Phi}{\rightarrow}}{\lim }F$.

\begin{lemma}
The functor of global co-sections is right exact.
\end{lemma}

\begin{pf}
Let $A$ be a co-sheaf on $\Phi $ with coefficients in $\cC $. The duality 
functor $\ ^0:\cC \to \cC ^0$ makes it into a sheaf $A^0$ with 
coefficients if $\cC ^0$. We have $A^{00}=A$ and 
$$(\underset{\underset{\Phi }{\rightarrow }}{\lim }A)^0=
\underset{\underset{\Phi }{\leftarrow }}{\lim }A^0.$$
Also the functor $A\mapsto A^0$ preserves exact sequences. 

Let 
$$0\to A \to B\to C\to 0$$
be an exact sequence of co-sheaves. Then 
$$0\leftarrow A^0\leftarrow B^0\leftarrow C^0\leftarrow 0$$
is also exact. Hence 
$$\underset{\underset{\Phi }{\leftarrow }}{\lim }A^0
\leftarrow 
\underset{\underset{\Phi }{\leftarrow }}{\lim }B^0
\leftarrow
\underset{\underset{\Phi }{\leftarrow }}{\lim }C^0
\leftarrow 0$$
is exact. Applying $\ ^0$ to the last sequence we obtain the desired exact 
sequence
$$\underset{\underset{\Phi }{\rightarrow }}{\lim }A\to 
\underset{\underset{\Phi }{\rightarrow }}{\lim }B\to 
\underset{\underset{\Phi }{\rightarrow }}{\lim }C\to 0.$$
\end{pf}

\begin{remark}
Assume that the category $\cC$ has enough projectives. Then so does
the category of cosheaves on $\Phi $. Indeed, let $W$ be a projective
object in $\cC$. Choose a cone $\sigma \in \Phi$ and denote by $J$
the "constant co-sheaf" $W$ on $\Star(\sigma )$ "extended by zero"
to $\Phi $. Then $J$ is a projective co-sheaf and every co-sheaf
is a quotient of a direct sum of such.
\end{remark}

Based on the above lemma and remark we may define the $i$-th
homology of a co-sheaf $\cV $ as the $i$-th left derived functor of
global co-sections:

$$H_i(\Phi ;\cV ):= H^{-i}\LL \lim_{\underset{\Phi }{\rightarrow}} \cV .$$

\begin{example}
Let $J$ be a co-sheaf as in the last remark. Then
$H_0(\Phi ;J)=W$ and $H_i(\Phi ;J)=0$, for $i>0$.
\end{example}

\subsection{Co-sheaf of sections with compact support}
Fix an orientation of each cone $\sigma $ in $\Phi $. Let $F$ be a sheaf
on $\Phi $. We want to define the complex $\Gamma _c(F)$ of co-sheaves  of
sections of $F$ with compact support. Let $\sigma \in \Phi$. Denote
by $i_{[\sigma ]}:[\sigma ]\hookrightarrow \Phi $ the open embedding
of the irreducible open set $[\sigma ]$. Denote by $F_{[\sigma ]}$ the
extension by zero $i_{[\sigma ] !}i_{[\sigma ]}^*F$.
 Consider the cellular
complex $C^\bullet(F_{[\sigma ]})$ and put

$$\Gamma _c(F)_{\sigma }:=C^\bullet (F_{[\sigma ]})[n].$$
This is a complex which is concentrated in degrees $[-d(\sigma),0]$.
For $\tau \leq \sigma$ we have  the inclusion of complexes
$\Gamma _c(F)_{\tau }\hookrightarrow \Gamma _c(F)_{\sigma }$ which makes
$\Gamma _c(F)$ a complex of co-sheaves.

The functor $\Gamma _c(\cdot )$ is exact, thus it extends trivially to
the derived category of sheaves on $\Phi $. Let $\Sh (\Phi )$ as usual 
denote the category of sheaves of  $\bbR$-vector spaces on $\Phi $.

\begin{prop} 
Assume that the fan $\Phi $ is complete. 
The functors $\LL{\underset{\Phi}{\varinjlim}}\Gamma _c(\cdot ) $ and
$\R \Gamma (\Phi ;\cdot)[n]$
are naturally isomorphic (as functors from the bounded derived category
$D^b(\Sh(\Phi ))$ to $D^b(\bbR -vect)$).
\end{prop}

\begin{pf}
Fix a sheaf $F$ on $\Phi $. Then $\Gamma _c(F)$ is a complex which consists
of direct sums of projective co-sheaves $J$ as in Example 6.3 above.
Thus the natural morphism
$$\LL \varinjlim \Gamma _c(F) @>>>\varinjlim \Gamma _c(F)$$
is an isomorphism.
 On the other hand
$\varinjlim \Gamma _c(F)$  is just the shifted
cellular complex
$C^\bullet(F)[n]$. Since the fan $\Phi $ is complete $C^\bullet(F)=
\R \Gamma(\Phi ;F)$. This proves the proposition.
\end{pf}

\subsection{Duality in $D_c^b(A_{\Phi }-mod)$}
We are going to define a version of Borel-Moore-Verdier duality for sheaves of
$A_{\Phi }$-modules.  Let us first
fix some notation.

Define the dualizing $A$-module as follows
$$\omega =\omega_{A/\bbR}:=A\otimes \det V^*.$$
Thus $\omega$ is a free $A$-module of rank one generated in
degree $2\dim V$, i.e. $\omega \simeq A(-2n)$.
For a complex $M$ of $A$-modules denote by $M^*$ the complex
$\Hom _A^\bullet (M,\omega )$. In what follows we will omit the
subscript $A$ (resp. $A_{\Phi }$) when talking about morphisms of
$A$-modules (resp. $A_{\Phi }$-modules).

\begin{defn}
Fix $F\in D_c^b(A_{\Phi }-mod)$.
The category of $A_{\Phi }$-modules has enough projective objects.
 Let $P(F)\to F$ be
a projective resolution of $F$. (We may assume that $P(F)$ is a finite
complex if projectives $A_{\Phi }$-modules).
 Then in particular for each cone
$\sigma \in \Phi$ the stalk $P(F)_{\sigma }$ is a complex of free
$A$-modules.
Consider the complex
of co-sheaves of $A$-modules $\Gamma _c(P(F))$. Finally put
$$\bbD(F):=\Gamma _c(P(F))^*[n].$$

This is a complex of sheaves of $A$-modules, whose stalk at $\sigma $ is
the complex
$$\bbD(F)_{\sigma }=\Gamma _c(P(F))_{\sigma }^*[n]=
C^\bullet (F_{[\sigma ]})^*.
 $$
Thus if $P(F)$ is a single sheaf in degree zero, then $\bbD(F)_{\sigma }$
is a complex concentrated in degrees $[-n, -n+d(\sigma )]$.
\end{defn}

\begin{remark}
Note that the duality $\bbD $ is a local functor, i.e. it commutes with 
restriction of sheaves to subfans.
\end{remark}

\begin{prop}
Let $\Phi $ be a fan in $V$. There is a natural isomorphism of functors 
from $D_c^b(A_{\Phi }-mod)$ to $D^b(A-mod)$:
$$\R \Gamma (\Phi ;\bbD (\cdot ))\simeq 
\R \Hom (C^\bullet (\cdot ),\omega ).$$ 
\end{prop}

\begin{pf}
Let $P$ be a projective $A_{\Phi }$-module. Then $\bbD (P)$ is a complex 
of sheaves with the following graded components
$$\bbD (P)^{-n+i}=\underset{d(\sigma )=i}{\bigoplus}
(P_{\sigma }^*)_{\Star (\sigma )},$$
where $(P_{\sigma }^*)_{\Star (\sigma )}$ is the extension by zero to $\Phi $ 
of the constant sheaf on $\Star (\sigma )$ with stalk $\Hom (P_{\sigma },
\omega )$. In particular, such sheaves are injective (when considered 
as sheaves of vector spaces) and hence 
$$\R \Gamma (\Phi ;\bbD (P)^{-n+i})=
\Gamma (\Phi ;\bbD (P)^{-n+i})=
\underset{d(\sigma )=i}{\bigoplus }P_{\sigma }^*.$$
Thus the complex of global sections $\R \Gamma (\Phi ;\bbD (P))$ is just the 
dual of the cellular complex $C^\bullet (P)$:
$$\R \Gamma (\Phi ;\bbD (P))=C^\bullet (P)^*.$$

The last formula remains true if $P$ is a complex of projective 
$A_{\Phi }$-modules. 

Now fix $F\in D_c^b(A_{\Phi }-mod)$ and let $P(F)\to F$ be its 
projective resolution. Since the cellular complex is an exact functor 
we have the quasi-isomorphism 
$$C^\bullet (P(F))\simeq C^\bullet (F).$$
Note that the complex $C^\bullet (P(F))$ consists of free $A$-modules, 
so it is a projective resolution of $C^\bullet (F)$. Thus
$$\R \Hom (C^\bullet (F), \omega )=\Hom (C^\bullet (P(F)), \omega )=
C^\bullet (P(F))^*.$$
Combining this with the equality 
$$\R \Gamma (\Phi ;\bbD (F))=\R \Gamma (\Phi ;\bbD (P(F)))=
C^\bullet (P(F))^*$$
we obtain the desired isomorphism.
\end{pf}

\begin{cor}
Let $\Phi $ be a fan in $V$. Assume that $\Phi $ is complete. 
Then there is a natural isomorphism of
functors from $D^b_c(A_{\Phi }-mod)$ to $D^b(A-mod)$:

$$\R\Gamma (\Phi ;\bbD (\cdot ))\simeq \R\Hom (\R \Gamma (\Phi ;\cdot ),
\omega
 ).$$
\end{cor}

\begin{pf}
This follows from the last proposition and Proposition 3.3
\end{pf}

Later on we will need the following version of the above corollary.

\begin{cor}
Let $\Phi $ be a fan in $V$ and 
 $\sigma \in \Phi $ be a cone. Consider the  subfan
 $\Psi :=[\partial \sigma ]\subset \Phi $.
Then there is a natural isomorphism of functors
from $D^b_c(A_{\Phi }-mod)$ to $D^b(A-mod)$:

$$\R\Gamma (\Psi ;\bbD (\cdot ))\simeq \R\Hom (\R \Gamma (\Psi ;\cdot ),
\omega
 )[n-d(\sigma )+1].$$
\end{cor}

\begin{pf} 
This follows from Proposition 6.7, Remark 6.6 and Corollary 3.5.
 \end{pf}

We are going to give an alternative description of the duality functor
(as in the usual case of sheaves on locally compact
spaces).

\begin{defn}
Define the {\it dualizing complex} $D_{\Phi }$ on the fan $\Phi $ as
$$D_{\Phi }:=\bbD (A_{\Phi }).$$
\end{defn}

The complex $D_{\Phi }$ is concentrated in degrees $[-n,o]$ and
$D_{\Phi }^{-n+i}$ is the direct sum over all cones $\tau $ of dimension $i$
of constant sheaves with stalk $A_{\tau }^*= \omega $
 on
the closed set $\Star (\tau )$ extended by zero to $\Phi $.

\begin{prop}
There is a canonical isomorphism of contravariant 
endofunctors of $D_c^b(A_{\Phi }-mod)$:
$$\bbD (\cdot )=\R \cH om (\cdot ,D_{\Phi }).$$
\end{prop}

\begin{pf}
Suppose that $F\in D_c^b(A_{\Phi }-mod)$ is a complex of projective 
$A_{\Phi }$-modules. 
 We are
going to show that complexes of $A_{\Phi }$-modules $\bbD (F)$ and
$\cH om (F,D_{\Phi })$ are naturally isomorphic.
To simplify the notation assume that $F$ is a single sheaf in degree 0.
Fix $\sigma \in \Phi$. Then
$$\begin{array}{rcl}
\cH om^{i-n}(F,D_{\Phi })_{\sigma } & = &
\Hom ^{i-n}(F_{[\sigma ]},D_{\Phi })\\
 & = & \Hom (F_{[\sigma ]},D_{\Phi }^{i-n})\\
 & = & \underset{\underset{\tau \leq \sigma }{d(\tau )=i}}{\bigoplus}
F_{\tau }^*\\
 & = &
\bbD (F)^{i-n}_{\sigma }.\\
\end{array}$$
This identification is compatible with the differential in the two
complexes. This proves the proposition.
\end{pf}

In fact the dualizing complex $D_{\Phi }$ is very simple. Namely, consider
$\omega $ as a constant sheaf on the open point $\Or \in \Phi$. Let
$\omega _{\Or }$ denote its extension by zero to $\Phi $. Note that the
stalk $D_{\Phi,\Or }$ is isomorphic to $\omega _{\Or }[n]$.

\begin{lemma}
The natural map
$$\omega _{\Or }[n]\rightarrow D_{\Phi }$$
is a quasi-isomorphism.
\end{lemma}

\begin{pf} The lemma claims that for each $\sigma \neq \Or $ the complex
$D_{\Phi ,\sigma }$ is acyclic. This complex is $C^\bullet (A_{[\sigma ]})^*$
and $C^\bullet (A_{[\sigma ]})$ is acyclic by Lemma 3.7.
\end{pf}

\subsection{Biduality}
Our next goal is to prove Theorem 6.23 below. Let us begin with some preparations.

\begin{defn}
A nonempty open subset $U\subset \Phi $ is {\em saturated} if whenever
the boundary of a cone is in $U$, then the cone itself is in $U$.
\end{defn}

\begin{example}
$\Phi $ and $\{ \Or \}$ are saturated.
\end{example}

\begin{defn}
For a nonempty open subset $U\subset \Phi $ define its {\em opposite}
$U^\prime $ as follows:
$$U^\prime :=\{ \sigma \in \Phi \vert \forall \tau \in U,\ \
\tau \cap \sigma =\Or \}.$$
\end{defn}

\begin{remark}
1.  $\Phi ^\prime =\{ \Or \}$, $\{ \Or \}^\prime =\Phi $.

2.  For any $U$ its opposite $U^\prime $ is saturated.

3.  We have $U\subset U^{\prime \prime }$.
\end{remark}

\begin{lemma}
If $U$ is saturated then $U=U^{\prime \prime }$.
\end{lemma}

\begin{pf}
Let $\Or \neq \sigma \in U^{\prime \prime }$. Let $\tau \leq \sigma $
be a face of dimension 1. Then $\tau \notin U^\prime $. Hence $\tau \in U$.
Thus by induction on dimension all faces of $\sigma $ are in $U$ and so
$\sigma $ is in $U$.
\end{pf}

\begin{cor}
The map $U\mapsto U^\prime $ is an involution of the collection of
saturated open subsets of $\Phi $.
\end{cor}

For an open subset $U\subset \Phi $ and an $A$-module $M$ denote as usual by
$M_U$ the extension by zero to $\Phi $ of the constant sheaf $M$ on $U$.
In case $U=\{ \Or \}$ we will also denote this sheaf by $M_{\Or }$.

\begin{remark}

Let $U\subset \Phi $ be open. Note the equality of
sheaves
$$\cH om(A_U,\omega _{\Or })=\omega _{U^\prime }.$$
Hence if $U$ is saturated, then by Lemma 6.17 the obvious map
$$A_U\rightarrow \cH om(\cH om (A_U, \omega _{\Or }), \omega _{\Or }),$$
$$a\mapsto (f\mapsto f(a))$$
is an isomorphism.
\end{remark}

\begin{defn}
Let us introduce the apropriate version of the cocellular complex
(see 3.14) for $A_{\Phi }$-modules. For a cone $\sigma \in \Phi $ let
$i_{\sigma }:\{ \sigma \}\hookrightarrow \Phi $ denote the inclusion.
Consider the dualizing module $A$-module $\omega $ as a sheaf on the
 point $\sigma $.
Then the
 $A_{\Phi }$-module $i_{\sigma *}\omega $ is a constant sheaf on $\Star
(\sigma )$ with stalk $\omega $. If $\tau \leq \sigma $ then there is
a natural surjection of sheaves $r_{\tau \sigma }:i_{\tau *}\omega \to
 i_{\sigma *}\omega $. Put
$$K_A^{-n+j}:=\underset{d(\sigma)=j}{\bigoplus } i_{\sigma *}\omega.$$
As usual, the maps $r_{\tau \sigma }$ with the sign $\pm 1$ define the
differential in the complex
$$K_A^\bullet :=K^{-n}\to K^{-n+1}\to ...\to K^0.$$
This complex is a resolution of the dualizing sheaf $\omega _{\Or }[n]$
(the proof is the same as that of Lemma 3.13).

Given an $A_{\Phi }$-module $F$ define its {\it cocellular complex}
$K_A^\bullet (F)$ as
$$K_A^\bullet (F):=\Hom ^\bullet (F,K_A^\bullet ).$$
\end{defn}

\begin{remark}
1. Note that $\Hom (F, i_{\sigma *}\omega)=\Hom (F_{\sigma }, \omega )$.
Hence if stalks $F_{\sigma }$ are free $A$-modules, then
$$K_A^\bullet(F)=\R \Hom (F, \omega _{\Or }[n]).$$

2. The cocellular complex $K_A^\bullet (F)$ is just the dual of
the cellular complex $C^\bullet (F)$:
$$K_A^\bullet (F)=\Hom ^\bullet (C^\bullet (F),\omega )=C^\bullet (F)^*.$$
\end{remark}

\begin{lemma}
Let $\sigma \in \Phi $ and put $U:=[\sigma ]^\prime$. Then
$$\cH om (A_U,\omega _{\Or })=\R \cH om(A_U,\omega _{\Or }).$$
\end{lemma}

\begin{pf}
Fix $\tau \in \Phi $. We need to show that the stalk $\R ^i\cH om(A_U,
\omega _{\Or })_{\tau }=0$ for $i>0$. This is the same as showing that
$\R ^i\Hom (A_W,\omega _{\Or })=0$, for $i>0$, where $W=U\cap [\tau ]$.
 This is trivial if $\tau =\Or $,
so assume that $\tau \neq \Or $.
By Remark 6.21(1)  above it suffices to show that the cocellular complex
$K_A^\bullet (A_W)$ is acyclic in  degrees $>-n$. If $W={\Or }$, then
$K_A^i (A_W)=0$ if $i>-n$. So assume that $W\neq {\Or }$. Then
the space $\tau -W$ is contractible. Therefore the cocellular
complex $K_A^\bullet (A_W)$ is isomorphic to the augmented (shifted)
 cochain complex of a
contractible space; hence it is acyclic.
\end{pf}

We are ready to prove the main result of this section.

\begin{thm}
Let $\Phi $ be a fan in $V$. 
The duality $\bbD$ is an anti-involution of the category
$D^b_c(A_{\Phi }-mod)$. More precisely, there exists an isomorphism
of functors
$$Id\overset{\sim}{\rightarrow}\bbD \circ \bbD .$$
\end{thm}

\begin{pf}
For $F\in D_c^b(A_{\Phi }-mod)$ we have a functorial morphism
$$\alpha : F\to \bbD \circ \bbD (F)=\R \cH om 
(\R \cH om (F,D_{\Phi }),D_{\Phi }
 ),$$
$$a\mapsto (f\mapsto f(a)).$$
Let $\sigma \in \Phi $. Then $A_{[\sigma ]}$ is a projective $A_{\Phi }$
-module. It suffices to prove that $\alpha (A_{[\sigma ]})$ is a
quasi-isomorphism for all $\sigma $.

By Remark 6.19 and Lemma 6.22 we have
$$\begin{array}{rcl}
 \bbD \circ \bbD (A_{[\sigma ]}) & = &
\R \cH om (\cH om (A_{[\sigma ]},
 \omega _{\Or }[n]), \omega _{\Or }[n])\\
 &  = &
\R \cH om(\omega _{[\sigma ]^\prime }, \omega _{\Or })\\
 &  = & \cH om(\omega _{[\sigma ]^\prime }, \omega _{\Or })\\
 &  = & A_{[\sigma ]}.\\
\end{array}$$

Hence $\alpha (A_{[\sigma ]})$ is an isomorphism by Remark 6.19. 
This proves the
theorem.
\end{pf}

\subsection{Poincare duality for $IH(\Phi )$}

\begin{prop}
Let $\Phi $ be a fan in $V$. 
The duality $\bbD $ preserves the subcategory ${\frak M} (\cA _{\Phi })$.
\end{prop}

\begin{pf}
Fix $\cM \in \cM (\cA _{\Phi })$. We need to show

a) $\bbD (\cM )$ is quasi-isomorphic to a single sheaf (sitting in
degree 0).

b) The $A _{\Phi }$-module $\bbD (\cM )$ is actually an $\cA _{\Phi }$-module
and as such is locally free.

c) The sheaf $\bbD (\cM )$ is flabby.

a), b):
Since the sheaf $\cM $ is flabby the cohomology of the complex
$\Gamma _c (\cM _{[\sigma ]})$
vanishes except in degree $-d(\sigma )$ (see Lemma 3.9).
Put $K:=H^{-d(\sigma )}\Gamma _c (\cM _{[\sigma ]})$; it
is a free $\cA_{\Phi,\sigma}$-module.
This means that $\R^i\Hom_A(\Gamma _c (\cM _{[\sigma ]}),\omega [n]) = 0$ for
$i\neq 0$ and $\R^0\Hom_A(\Gamma _c (\cM _{[\sigma ]}),\omega [n])=\Ext
^i_A(K,\omega )$ is a free
$\cA_{\Phi,\sigma}$-module. This proves a) and b).

It remains to prove
that the restriction map
$\bbD (\cM )_{\sigma }\to \Gamma (\partial \sigma ,\bbD (\cM ))$ is
surjective (Lemma 3.1). Consider the open set $[\partial \sigma ]$ and the
cellular complex
$C^\bullet (\cM _{[\partial \sigma ]})$. By Corollary 3.5 
it is quasi-isomorphic to
$\R \Gamma (\partial \sigma ;\cM )[-n+d(\sigma )-1]$.
Since $\cM $ is flabby we have $H^j(\partial \sigma ;\cM)=0$
unless $j=0$. Put $N:= H^0
(\partial \sigma ;\cM)$. Since $\cM $ is locally free we
have $\Ext ^j_A(N,A)=0$ unless $j=n-d(\sigma )+1$.
Moreover by Corollary 6.9
$$\Gamma (\partial \sigma ; \bbD (\cM ))=\R \Hom (\R \Gamma
(\partial \sigma ;\cM), \omega)[n-d(\sigma )+1].$$

Thus $\Gamma (\partial \sigma ; \bbD (\cM ))=\Ext^{n-d(\sigma )+1}(N,\omega)$.

We have the exact sequence of $A$-modules

$$0\to K \to \cM _{\sigma } \to N \to 0,$$
which gives rise to the exact sequence

$$0\to \Ext _A^{n-d(\sigma )}(\cM _{\sigma }, \omega )\to
 \Ext _A^{n-d(\sigma )}(K,\omega )
\overset{f}{\rightarrow } \Ext _A^{n-d(\sigma )+1}(N,\omega ) \to 0.$$
The map $f$ coincides with the restriction map

$$\bbD (\cM )_{\sigma }\to \Gamma (\partial \sigma ;\bbD (\cM )),$$
hence the latter is surjective. This proves c).

\end{pf}

\begin{cor}
Let $\Phi $ be a fan in $\bbR ^n$.
The duality functor induces an anti-involution of the category
${\frak M}(\cA _\Phi )$.
\end{cor}

\begin{pf} Since $\bbD $ is an anti-involution of $D_c^b(A_{\Phi }-mod)$ and
it preserves ${\frak M}(\cA _{\Phi })$ the corollary follows.
 \end{pf}

\begin{cor}
Let $\Phi $ be a fan in $\bbR ^n$ and $\sigma \in \Phi $.
 Let $\cL ^{\sigma }_{\Phi }(k)$ be the minimal sheaf
on $\Phi $ based at $\sigma $. Then $\bbD (\cL ^{\sigma }_{\Phi }(k))\simeq
\cL ^{\sigma }_{\Phi }(-k-2d(\sigma ))$. In particular, 
$\bbD (\cL _{\Phi })\simeq \cL _{\Phi }$.
\end{cor}

\begin{pf}
By the above corollary the duality sends indecomposable objects of
${\frak M}(\cA _\Phi )$ to indecomposable ones. We have
$\supp (\bbD (\cL ^{\sigma }_{\Phi }(k)))\subset \Star (\sigma )$.
Also $\bbD (\cL ^{\sigma }_{\Phi }(k))_{\sigma }=
 \Ext _A^{n-d(\sigma )}(\cA _{\Phi ,\sigma }(k),\omega )
\simeq \cA _{\Phi ,\sigma }(-k-2d(\sigma ))$ (see the proof of 
Proposition 6.24 above).
 Thus the corollary follows from Proposition 5.2(3) and Theorem 5.3.
\end{pf}

\begin{cor}
Let $\Phi $ be a complete fan in $\bbR ^n$. Then
there exists an isomorphism of $A$-modules:
$$\Gamma (\Phi ;\cL _{\Phi })\simeq \Hom (\Gamma (\Phi ;\cL _{\Phi }),\omega ),$$
i.e. the free $A$-module $\Gamma (\Phi ;\cL _{\Phi })$ is self-dual.
 \end{cor}

\begin{pf}
The minimal sheaf $\cL _{\Phi }$ is flabby, so $H^i(\Phi ;\cL _{\Phi })=0$
for $i>0$. By Theorem 4.7 the $A$-module $\Gamma (\Phi ;\cL _{\Phi })$ 
is free. The
above corollary implies that the same is true for $\bbD (\cL _{\Phi })$.
Hence the natural isomorphism of Proposition 6.8
$$\R \Gamma (\Phi ;\bbD (\cL _{\Phi }))=\R \Hom (\R \Gamma (\Phi ;
\cL _{\Phi }), \omega )$$
reduces to
$$\Gamma (\Phi ;\bbD (\cL _{\Phi }))=\Hom (\Gamma (\Phi ; \cL _{\Phi }),
\omega ).$$

A choice of an isomorphism of $\bbD (\cL _{\Phi })$ and $\cL _{\Phi }$
provides an isomorphism of $A$-modules
$$\Gamma (\Phi ;\cL _{\Phi })\simeq \Hom (\Gamma (\Phi ;\cL _{\Phi }),
\omega).$$
 \end{pf}

Consider the 1-dimensional (graded) vector space $\overline{\omega }$. It has
degree $2n$.

\begin{cor}
Let $\Phi $ be a complete fan in $\bbR ^n$. There exists an isomorphism
of graded vector spaces

$$IH(\Phi )\simeq \Hom _{{\Bbb R}}(IH(\Phi ),\overline {\omega }).$$
\end{cor}

\begin{pf}  Immediate from the previous corollary.
\end{pf}

\begin{cor}
Let $\Phi $ be a complete fan in ${\Bbb R}^n$. Then

1. $ih_{n-j}(\Phi )=ih_{n+j}(\Phi )$ for all $j$. 

2. $ih_j(\Phi )=0$ unless $j$ is even and $j\in [0,2n]$.

3. $ih_0(\Phi )=1=ih_{2n}(\Phi )$.
\end{cor}

\begin{pf}
The first assertion is immediate from the last corollary. Recall that
$ih_0(\Phi )=1$ and $ih_j(\Phi )=0$ if $j<0$ or $j$  odd (Lemma 5.7). 
This implies
the last two assertions.
\end{pf}

\section{Toward Hard Lefschetz and the combinatorial invariance}
\label{section:HL}

Throughout this section $\Phi$ will denote a complete fan in a vector
space $V$ of dimension $n$.

\subsection{Ampleness in the context of fans}
Consider the short exact sequence of sheaves
\[
0 @>>> \Omega^1_\Phi @>>> V^\ast @>>> \cG @>>> 0\ ,
\]
where $V^\ast$ denotes the constant sheaf and $\Omega^1_{\Phi,\sigma} =
\Span(\sigma)^\perp$. Since constant sheaves have trivial higher cohomology and
$\Omega^1_\Phi$ is supported on $\Phi_{\leq n-1}$, the long exact sequence
in cohomology reduces, in low degrees, to the short exact sequence of vector
spaces
\[
0 @>>> V^\ast @>>> \Gamma(\Phi,\cG) @>>> H^1(\Phi;\Omega^1_\Phi) @>>> 0\ .
\]
The space $\Gamma(\Phi;\cG)$ consists of continuous, cone-wise linear functions
on $\Phi$.

For any object $\cM$ of ${\frak M}(\cA_\Phi)$, the elements of
$\Gamma(\Phi;\cG)$ act naturally on the free graded $A$-module
$\Gamma(\Phi;\cM)$ by endomorphisms of degree two. Clearly, the
induced action on the graded vector space $\overline{\Gamma(\Phi;\cM)}$
factors through $H^1(\Phi;\Omega^1_\Phi)$.

\begin{defn}
An element $\overline l$ of $H^1(\Phi;\Omega^1_\Phi)$ is called {\em ample}
iff it admits a lifting $l\in\Gamma(\Phi;\cM)$ which is strictly convex.
\end{defn}

\subsection{Hard Lefschetz for complete fans}
The statement of Conjecture \ref{conj:HL} (below) is the 
analog of the Hard Lefschetz Theorem
in the present context. Recall that $\cL_\Phi$ denotes the indecomposable
 object of
 ${\frak M}(\cA_\Phi)$ based at the origin, with $\cL_{\Phi,\Or}=
\bbR$. For a graded vector space $W$ we will denote by $W^{(i)}$ the subspace
of homogeneous elements of degree $i$.

\begin{conj}[Hard Lefschetz Conjecture]\label{conj:HL}
An ample $\overline{l}\in H^1(\Phi;\Omega^1_\Phi)$ induces a Lefschetz
operator on the graded vector space $IH(\Phi)$, i.e. for every $i$ the map
$\overline{l}^i : IH(\Phi)^{(n-i)} @>>> IH(\Phi)^{(n+i)}$ is an isomorphism.
\end{conj}

For a rational fan $\Phi $ this conjecture follows immediately from 
results in [BL], ch.15.

The above conjecture has the following standard corollary.

\begin{cor}\label{cor:HL} Assume the Hard Lefschetz Conjecture. Then for 
an ample $\overline{l}$ 
the map $\overline{l}: IH(\Phi)^{(i)} @>>> IH(\Phi)^{(i+2)}$ is injective for
$i\leq n-1$ and surjective for $i\geq n-1$. In particular 
$$ih_0(\Phi )\leq ih_2(\Phi )\leq ...\leq ih_{2[n/2]}(\Phi ).$$
\end{cor}

\subsection{Local intersection cohomology.}
Let $\Psi $ be a fan in $V$ and $\sigma \in \Psi$. 

\begin{defn}
Define the {\it local intersection cohomology} space
$$IP(\sigma ):=\overline{\cL _{\Psi ,\sigma }}$$
and denote by $ip(\sigma )$ the corresponding Poincare polynomial.
\end{defn}

\begin{remark}
1. Consider the subfan $[\sigma ]\subset \Psi $. Then by Proposition 
5.2(4) 
$\cL _{\Psi }\vert _{[\sigma ]}=\cL _{[\sigma ]}.$ Thus in particular 
the local intersection cohomology $IP(\sigma )$ depends only on the cone 
$\sigma $ and not on the fan $\Psi $. 

2. Note that $ip_j(\sigma )=0$ if $j$ is odd or negative.
\end{remark}

\subsection{The global-local formula.}
For a cone $\sigma \subset V$ of dimension $d(\sigma )=d+1\geq 2$ 
consider the subspace $W:=\Span (\sigma )\subset V$.  
Choose a linear isomorphism 
$$W\simeq {\Bbb R}^{d}\times {\Bbb R}$$
so that the ray $(0,{\Bbb R}^+)$ lies in the interior of $\sigma $. 
This defines the projection $p$ of $\sigma $ to ${\Bbb R}^{d}$ so 
that the image $\overline{\partial \sigma }:=p(\partial \sigma )$ is a 
complete fan in ${\Bbb R}^{d}$. Also $\partial \sigma $ becomes the graph 
of a function $l:{\Bbb R}^{d} \to {\Bbb R}$ which is piecewise linear 
and strictly convex with respect to the fan $\overline{\partial \sigma }$. 
In particular $\overline{l}\in H^1(\overline{\partial \sigma }; \Omega ^1 
_{\overline{\partial \sigma }})$ is an ample class. 

Note that $A_W=A_{{\Bbb R}^{d}}[l]$. Let $m_1\subset A_W$ and $m_2\subset 
A_{{\Bbb R}^{d}}$ denote the maximal ideals, so that $m_1=m_2[l]$. 

The minimal sheaf $\cL _{[\partial \sigma ]}$ as the 
$\cA _{[\partial \sigma ]}$-module is obtained by extension of scalars
$$\cL _{[\partial \sigma ]}=\cA _{[\partial \sigma ]}\otimes 
_{\cA _{\overline{\partial \sigma }}}\cL _{\overline{\partial \sigma }}.$$
Thus the intersection cohomology $IH(\overline{\partial \sigma })$ is an 
${\Bbb R}[l]$-module. We have

$$IH(\overline{\partial \sigma })/l\cdot IH(\overline{\partial \sigma })
\simeq \Gamma (\partial \sigma ;\cL _{[\partial \sigma ]})/
m_1\Gamma (\partial \sigma ;\cL _{[\partial \sigma ]})
\simeq \overline{\cL _{[\sigma ],\sigma }}=IP(\sigma ).$$

The Hard Lefschetz Conjecture (for $\overline{\partial \sigma }$) 
implies that $\overline{l}$ is a Lefschetz operator on 
$IH(\overline{\partial \sigma })$. Thus $IP(\sigma )$ is 
isomorphic to the $l$-primitive part of $IH(\overline{\partial \sigma })$. 
In particular the 
Poincare polynomial $ih(\overline{\partial \sigma })$ depends only on 
$\sigma $ and not on a particular choice of the isomorphism
$$\Span (\sigma )\simeq {\Bbb R}^{d}\times {\Bbb R}.$$

Let us summarize our discussion in the following corollary.

\begin{cor}
Let $\sigma \subset V$ be a cone of dimension $d+1\geq 2$. 
Choose an isomorphism 
$$\Span (\sigma )\simeq {\Bbb R}^{d}\times {\Bbb R}$$
as above, so that the projection $\overline{\partial \sigma }$ of 
the fan $[\partial \sigma ]$ is a complete 
fan in ${\Bbb R}^{d}$. Then the Hard Lefschetz Conjecture implies that 

1. The Poincare polymonial $ih(\overline{\partial \sigma })$ 
is independent of the choices made.

2. $$ip_j(\sigma ) = \left\lbrace
\begin{array}{cc} 
ih_j(\overline{\partial \sigma}) - ih_{j-2}(\overline{\partial \sigma })	& \text{for $0\leq j\leq d$} \\
0				& \text{otherwise}
\end{array}
\right.$$
\end{cor}

\begin{defn}
In the above corollary denote the Poincare polynomial 
$ih(\overline{\partial \sigma })$ by $ih(\sigma )$. We have 
$$ip_j(\sigma ) = \left\lbrace
\begin{array}{cc} 
ih_j(\sigma ) - ih_{j-2}(\sigma )	& \text{for $0\leq j\leq d$} \\
0				& \text{otherwise}
\end{array}
\right.$$
We call the last equation the {\it global-local formula}.
In case $\sigma $ has dimension 0 or 1  put $ih(\sigma )=1$. 
(Note that in this case $ip(\sigma )=1$.)  
Note also that if cones $\sigma $ and $\sigma ^\prime$ are linearly 
isomorphic, 
then $ip(\sigma )=ip(\sigma ^\prime )$ and $ih(\sigma )=ih(\sigma ^\prime )$. 
\end{defn}

\subsection{The local-global formula}
In this section we express the Poincar\'e polynomial $ih(\Phi )$ of 
a complete fan $\Phi $ in terms of the local Poincar\'e polynomials 
$ip(\sigma )$ for $\sigma \in \Phi $. The argument is standard and is 
independent of any conjectures.

\begin{prop}\label{prop:ih-ip}
For a complete fan $\Phi$ in ${\Bbb R}^n$ we have the following relation 
between 
Poincar\'e polynomials in the variable $q$:
\[
ih(\Phi)(q) = \sum_{\sigma\in\Phi} (q^2-1)^{n-\dim\sigma }ip(\sigma)(q).
\]

\end{prop}
\begin{pf}
The quasi-isomorphism (Proposition 3.3) 
\[
\Gamma(\Phi;\cL_\Phi) @>>> C^\bullet(\cL_\Phi)
\]
induces the quasi-isomorphism
\[
\Gamma(\Phi;\cL_\Phi)\otimes^\bbL_A\bbR @>>> C^\bullet(\cL_\Phi)\otimes^\bbL_A\bbR
\]
and the equality of the graded Euler characteristics
\[
\chi(\Gamma(\Phi;\cL_\Phi)\otimes^\bbL_A\bbR) =
\chi(C^\bullet(\cL_\Phi)\otimes^\bbL_A\bbR)\ .
\]

Since $\Gamma(\Phi;\cL_\Phi)$ is free over $A$, the canonical map
\[
\Gamma(\Phi;\cL_\Phi)\otimes^\bbL_A\bbR @>>> \Gamma(\Phi;\cL_\Phi)\otimes_A\bbR
= IH(\Phi)
\]
is a quasi-isomorphism and
\[
\chi(\Gamma(\Phi;\cL_\Phi)\otimes^\bbL_A\bbR) = ih(\Phi)\ .
\] 

Since $C^\bullet(\cL_\Phi))$ is a complex of finitely generated $A$-modules
and $A$ has finite $\Tor$-dimension it follows that
\[
\chi(C^\bullet(\cL_\Phi)\otimes^\bbL_A\bbR ) =
\sum_i (-1)^i\chi(C^i(\cL_\Phi)\otimes^\bbL_A\bbR)
\]
Since $C^i(\cL_\Phi)$ is isomorphic to
$\bigoplus_{\dim\sigma = n-i}\cL_{\Phi,\sigma}$ the above formulas imply the
equality
\[
ih(\Phi) = \sum_{\sigma\in\Phi}(-1)^{n-\dim\sigma}
\chi(\cL_{\Phi,\sigma}\otimes^\bbL_A\bbR)\ .
\]
By definition of $\cL_\Phi$, the stalk $\cL_{\Phi,\sigma}$ is a free module
over $\cA_{\Phi,\sigma}$ of graded rank $ip(\sigma)$. The standard caclulation
with the Koszul complex shows that $\cA_{\Phi,\sigma}\otimes^\bbL_A\bbR$ is
represented by the complex (with trivial differential)
$\bigoplus_i{\bigwedge}^i\sigma^\perp$.
It follows that
\[
\chi(\cL_{\Phi,\sigma}\otimes^\bbL_A\bbR) = (1-q^2)^{n-\dim\sigma}ip(\sigma)
\]
and
\[
ih(\Phi) = \sum_{\sigma\in\Phi}(-1)^{n-\dim\sigma}
(1-q^2)^{n-\dim\sigma}ip(\sigma) =
\sum_{\sigma\in\Phi} (q^2-1)^{n-\dim\sigma }ip(\sigma)
\]
\end{pf}

\subsection{Summary.}
Assuming the Hard Lefschetz Conjecture for complete fans we have associated 
two polynomials $ip(\sigma )$ and $ih(\sigma )$ to any cone 
$\sigma \subset V$. The odd coefficients of these polynomials vanish and the 
following relations hold (here $d+1:=d(\sigma )$):

1. $ip(\Or )=ih(\Or )=1.$

2. $ip_j(\sigma ) = \left\lbrace
\begin{array}{cc} 
ih_j(\sigma ) - ih_{j-2}(\sigma )	& \text{for $0\leq j\leq d$} \\
0				& \text{otherwise}
\end{array}
\right.$
 
3. $ih(\sigma )(q) = \underset{\tau <\sigma}{\sum} 
(q^2-1)^{d-d(\tau) }ip(\tau)(q).$

Indeed, the first two relations are contained in Definition 7.7
 and the third one 
follows from Proposition 7.8 applied to the complete fan $\overline{\partial 
\sigma }$ as in Corollary 7.6 above. 

As an immediate consequence of the above relations we obtain (by induction 
on the dimension $d(\sigma )$) that the polynomials $ip(\sigma )$ and 
$ih(\sigma )$ are combinatorial invariants of $\sigma $, i.e. they depend 
only on the face lattice of $\sigma $. 

Recall that in case $d> 0$ the polynomial $ih(\sigma )$ is defined 
as $ih(\overline{\partial \sigma })$ for a complete fan 
$\overline{\partial \sigma }$ of dimension $d$. Hence it follows from 
Corollary 6.29  
and Corollary 7.3 that

1. $ih_0(\sigma )=1=ih_{2d}(\sigma ).$

2. $ih_j(\sigma )=0$, unless $j$ is even and $j\in [0,2d]$. 

3. For all $j$ $ih_{d-j}(\sigma )=ih_{d+j}(\sigma )$. 

4. $ih_0(\sigma )\leq ih_2(\sigma )\leq ...\leq ih_{2[d/2]}(\sigma )$.

\subsection{The $h$-vector and Stanley's conjectures.}
Let $Q\subset {\Bbb R}^n$ be a convex polytope of dimension $d$. In [S]  
Stanley defined two polynomials $g(Q)$ and $h(Q)$. These polynomials 
are defined simultaneously and recursively for faces of $Q$, including the 
empty face $\emptyset$, as follows:

1. $g(\emptyset )=h(\emptyset )=1.$

2. $g_j(Q) = \left\lbrace
\begin{array}{cc} 
h_j(Q) - h_{j-1}(Q)	& \text{for $0\leq j\leq [d/2]$} \\
0				& \text{otherwise}
\end{array}
\right.$

3. $h(Q)(t)=\sum_{P\subset Q}(t-1)^{d-d(P)-1}g(P)(t),$
where the last summation is over all proper faces $P$ of $Q$ including 
the empty face $\emptyset $. Here $d(P)$ is the dimension of $P$ and 
$d(\emptyset )=-1$.

Stanley proved (in a more general context of Eulerian posets) the 
``Poincar\'e duality'' for $h(Q)$:
$$h_j=h_{d-j},$$
and conjectured that 
$$0\leq h_0\leq h_1 \leq ...\leq h_{[d/2]}.$$

Let us show how this conjecture follows from the Hard Lefschetz Conjecture. 
Namely, consider the space ${\Bbb R}^n$ (which contains $Q$) as a hyperplane 
$({\Bbb R}^n,1)\subset {\Bbb R}^{n+1}$. Let $\sigma \subset {\Bbb R}^{n+1}$ 
be the cone with vertex at the origin $\Or$ which is spanned by $Q$. Then 
$d(\sigma )=d+1$. Nonempty faces of $\sigma $ are in bijective 
correspondence with 
faces of $Q$ (with a shift of dimension by 1), where the origin $\Or $ 
corresponds to the empty face $\emptyset \subset Q$.
Assuming the Hard Lefschetz Conjecture the polynomials $ih(\sigma )$ and 
$ip(\sigma )$ are defined, and by induction on dimension one concludes that 
$$ih(\sigma )(q)=h(Q)(q^2),\ \ \ ip(\sigma )(q)=g(Q)(q^2).$$
Thus Stanley's conjecture follows from the corresponding statement about the 
coefficients of $ih(\sigma )$.

\section{Calai conjecture (after T.~Braden and R.~MacPherson)}
\label{section:Kalai}
The statement of following theorem is the $ip$-analog of the
inequalities conjectured by G.~Kalai and proven, in the rational case,
by T.~Braden and R.D.~MacPherson in \cite{BM}. Our proof follows the
same pattern as the one in \cite{BM}. However, major simpifications
result from abscence of rationality hypotheses and, consequently,
any ties to geometry whatsoever. 

Suppose that $\sigma$ is a cone in $V$ and let $[\sigma]$ denote as usual the
corresponding ``affine'' fan which consists of $\sigma$ and all of it's
faces. Let $\tau \leq \sigma $ be a face. We know (5.2(4)) that 
$\cL _{[\sigma ]}\vert _{[\tau ]} = \cL _{[\tau ]}.$ Recall the graded 
vector spaces (7.4)
$$IP(\sigma )=\overline{\cL _{[\sigma ],\sigma }}, \ \ \ IP(\tau )
=\overline{\cL _{[\tau ],\tau }}$$
and the corresponding Poincar\'e polynomials $ip(\sigma )$, $ip(\tau )$. 
Consider the minimal sheaf $\cL _{[\sigma ]}^{\tau }\in 
{\frak M}(\cA _{[\sigma ]})$. Its support is $\Star (\tau )$ and we put 
$$IP(\Star (\tau )):=\overline{\cL ^{\tau }_{[\sigma ], \sigma }}.$$
Let $ip(\Star (\tau ))$ denote the corresponding Poincar\'e polynomial.

\begin{thm}\label{thm:Kalai}
Suppose that $\sigma$ is a cone (in $V$) and $\tau$ is a face of $\sigma$.
Then, there is an inequality, coefficient by coefficient, of polynomials with
non-negative coefficients
\[
ip(\sigma)\geq ip(\tau)\cdot ip(Star(\tau))\ .
\]
\end{thm}

\begin{pf} Let
$\iota : \Star(\tau) @>>> \Phi$ denote the closed embedding.
Then $\iota _*\iota ^{-1}\cL _{[\sigma ]}\in {\frak M} (\cA _{[\sigma ]})$. 
Indeed, the sheaf $\iota ^{-1}\cL _{[\sigma ]}$ is flabby, hence so is 
$\iota _*\iota ^{-1}\cL _{[\sigma ]}$. Moreover, $\iota _*$ is the extension 
by zero, so $\iota _*\iota ^{-1}\cL _{[\sigma ]}$ is locally free.

Thus by the structure Theorem 5.3 there is a direct sum decomposition
\[
\iota _*\iota ^{-1}\cL_{[\sigma]}
\simeq \bigoplus_{\rho\geq\tau}\cL_{[\sigma]}^\rho\otimes V_\rho
\]
where the multiplicities $V_\rho$ are certain graded vector spaces.
Comparing the stalks at $\tau $ we find that 
\[
\cL_{[\sigma ],\tau }\simeq \cL _{[\sigma ],\tau }^{\tau }\otimes V_{\tau }.
\]
Hence $V_{\tau }=IP(\tau )$. 

On the other hand, comparing the stalks at $\sigma $ we find
\[
\cL_{[\sigma],\sigma}
\simeq
\cL_{[\sigma],\sigma }^\tau \otimes V_\tau
\oplus
\bigoplus_{\rho > \tau}\cL_{[\sigma], \sigma }^\rho \otimes V_\rho\ .
\]
In particular
\[
IP(\sigma )
\simeq 
IP(\Star (\tau ))\otimes IP(\tau )
\oplus
\bigoplus_{\rho > \tau}\overline{\cL_{[\sigma], \sigma }^\rho}
\otimes V_\rho\ .
\]

Numerically this amounts to the inequality
\[
ip(\sigma)\geq ip(\tau)ip(\Star(\tau))\ .
\]
\end{pf}

\end{document}